\input amstex
\documentstyle{amsppt}
\magnification1200
\tolerance=10000
\overfullrule=0pt
\def\n#1{\Bbb #1}
\def\p{\Bbb C_{\infty}}

\def\min{\hbox{min }}

\def\ord{\hbox{ord }}

\def\De{\Delta}

\def\e11{E_{11}}

\def\ga{\goth A}

\def\ve{\varepsilon}
\def\de{\delta}

\def\ga{\gamma}

\def\be{\beta}
\def\th{\theta}
\def\al{\alpha}

\def\g{\goth }

\topmatter
\title
Calculation of $h^1$ of some Anderson t-motives
\endtitle
\author
S. Ehbauer, A. Grishkov, D. Logachev\footnotemark \footnotetext{E-mail: logachev94{\@}gmail.com \phantom{*********************************************************}}
\endauthor
\NoRunningHeads
\address
First and third authors: Departamento de Matem\'atica, Universidade Federal do Amazonas, Manaus, Brasil.
Second author: Departamento de Matem\'atica e estatistica
Universidade de S\~ao Paulo. Rua de Mat\~ao 1010, CEP 05508-090, S\~ao Paulo, Brasil, and Omsk State University n.a. F.M.Dostoevskii. Pr. Mira 55-A, Omsk 644077, Russia.
\endaddress
\abstract We consider Anderson t-motives $M$ of dimension 2 and rank 4 defined by some simple explicit equations parameterized by $2\times2$ matrices. We use methods of explicit calculation of $h^1(M)$ --- the dimension of their cohomology group $H^1(M)$ ( = the dimension of the lattice of their dual t-motive $M'$) developed in our earlier paper. We calculate $h^1(M)$ for $M$ defined by all matrices of the form $\left(\matrix 0& a_{12} \\  a_{21}& 0 \endmatrix \right)$, and by some matrices of the form $\left(\matrix a_{11}& a_{12} \\  a_{21}& 0 \endmatrix \right)$. These methods permit to make analogous calculations for most (probably all) t-motives.

$h^1$ of all Anderson t-motives $M$ under consideration satisfy the inequality $h^1(M)\le4$, while in all known examples we have $h^1(M)=0,1,4$. Do exist $M$ of this type having $h^1=2,3$? We do not know, this is a subject of further research.
\endabstract
\endtopmatter
\document
{\bf 0. General introduction.} Anderson t-motives ([G], 5.4.2, 5.4.18, 5.4.16) are the function field analogs of abelian varieties (more exactly, of abelian varieties with multiplication by an imaginary quadratic field (of MIQF-type), see for example [L]). Nevertheless, this analogy is not complete. For example, let $M$ be an Anderson t-motive, resp. $A$ an abelian variety. We can associate them a lattice $L(M)$, resp. $L(A)$. For abelian varieties, the functor $A\mapsto L(A)$ has a good description, see below. For Anderson t-motives, the situation is much worse, for example, $h_1(M)$ --- the dimension of $L(M)$ --- can be less than it is expected to be. We started a study of the lattice map of Anderson t-motives in [GL17], [GL18]. We proved in [GL17] that (roughly speaking) in a system of neighborhoods of a fixed Anderson t-motive the lattice map is an isomorphism. We developed in [GL18] a method of calculation of $h_1(M)$, as well as of $h^1(M)$ --- the dimension of the cohomology group of $M$, and we gave an example that not always $h_1(M)=h^1(M)$ (unlike the case of abelian varieties).
\medskip
The present paper is a continuation of [GL18]. We apply the method of calculation of $h^1(M)$ of [GL18] to a larger class of Anderson t-motives. We show that this method can be applied for most (probably all) t-motives.
\medskip
The cases considered in the present paper form a tiny part of the whole problem of finding of $h^1(M)$ for all $M$. This whole problem is really enormous. Clearly it cannot be solved without use of computers. Unfortunately writing of a corresponding program is not an easy task, because there exists a large diversity of cases. It is hardly likely that one scientist will be able to solve the problem, it requires a work of a large team of scientists. The authors hope that the results of the present paper will stimulate further research.
\medskip
The structure of the paper is the following. Sections 1 --- 4 are introductory. Section 2 contains definitions on Anderson t-motives, Section 3 gives necessary results of [GL18], and Section 4 explains general methods of calculation. For the conjectures, problems of further research and justification of the subject see Theorem 1.6; Conjecture 1.8; Problem 1.12 and the end of Section 1; Section 3.6; Conjecture 3.7; Section 4.7. In Sections 5, 6 we consider the case of $A$ of the form (3.4) (all cases), and in Sections 7, 8 we consider the case of $A$ of the form (3.5) (some cases). In Section 9 we give some calculations for the cases of $A$ of the form (3.5) which are not considered in Sections 7, 8. They can be useful for future researchers.
\medskip
{\bf 1. More detailed introduction.} Here we give more details. First, we recall the number field case. For an abelian variety $A$ of dimension $g$ there exist its homology and cohomology groups $H_1(A, \n Z)$ and $H^1(A, \n Z)$ (both these groups are isomorphic to $\n Z^{2g}$), and a $\n Z$-perfect pairing between them:

$$H_1(A, \n Z)\otimes_{\n Z}H^1(A, \n Z)\to \n Z\eqno{(1.1)}$$

There is an inclusion $\ga: H_1(A, \n Z)\to \n C^g$ such that $H_1(A, \n Z)$ forms a lattice in $\n C^g$. We have $A=\n C^g/H_1(A, \n Z)$. We have
\medskip
{\bf Theorem 1.2.} Abelian varieties of dimension $g$ over $\n C$ are in 1 -- 1 correspondence with $\n Z$-lattices of dimension $2g$ in $\n C^g$, satisfying the Riemann condition.
\medskip
Finally, for an abelian variety $A$ we can define its dual variety $A'$. There exist canonical isomorphisms

$$H_1(A, \n Z)\to H^1(A', \n Z), \ \ H^1(A, \n Z)\to H_1(A', \n Z)\eqno{(1.3)}$$

Let us give necessary definitions for the case of Anderson t-motives. Let $q$ be a power of a prime $p$, $\n F_q$ the finite field of order $q$. The function field analog of $\n Z$ is the ring of polynomials $\n F_q[\th]$ where $\th$ is an abstract variable. The analog of the archimedean valuation on $\n Q$ is the valuation at infinity on the fraction field $\n F_q(\th)$ of $\n F_q[\th]$; it is denoted by $ord$, it is uniquely determined by the property $\ord(\th)=-1$. The completion of an algebraic closure of the completion of $\n F_q(\th)$ with respect the valuation "ord" is the function field analog of $\n C$. It is denoted by $\p$.
\medskip
Abelian varieties have one discrete invariant --- their dimension $g$. Unlike them, Anderson t-motives have two invariants: dimension and rank (see 1.2.2); Anderson t-motives of dimension $n$ and rank $r$ are analogs of abelian varieties of dimension $r$ of MIQF-type, of signature $(n,r-n)$.
\medskip
An Anderson t-motive $M$ has the homology and cohomology groups (see [G], 5.9.11 (2), (3)\footnotemark \footnotetext{Goss uses a notation $H_1(E)$ instead of $H_1(M)$. This is practically the same: there is a 1 -- 1 correspondence between t-modules $E$ and t-motives $M$.} and Definition 2.5 of the present paper) $H_1(M, \n F_q[T])=H_1(M)$ and $H^1(M, \n F_q[T])=H^1(M)$ which are free $\n F_q[T]$-modules (here $T$ is an abstract variable, it is one of the generators of the Anderson ring, see Definition 2.1). The ranks of $H_1(M)$, $H^1(M)$ are denoted by $h_1(M), \ h^1(M)$ respectively. By analogy with the number field case we can expect that always $h_1(M)=h^1(M)=r$. But unlike the case of abelian varieties, they can be less than $r$.
\medskip
Like for the case of abelian varieties, for an Anderson t-motive $M$ there exists the dual t-motive $M'$ (see [GL07]). Analogs of (1.3) hold for Anderson t-motives (see, for example, [GL18], Proposition 1.9): there exist canonical isomorphisms

$$H_1(M)\to H^1(M'), \ \ H^1(M)\to H_1(M')$$

In particular, a method of calculation of $h^1$ permits us to calculate the $h_1$ as well: we apply it to the dual t-motive.

Analog of (1.1) is a pairing

$$\pi: H_1(M)\otimes_{\n F_q[T]} H^1(M)\to \n F_q[T]\eqno{(1.4)}$$

Counterexample of [GL18] shows that not always $h_1(M) = h^1(M)$, hence (1.4) is not always perfect. There is
\medskip
{\bf Theorem 1.5.} (Anderson, [A]; [G], 5.9.14). $h^1(M)=r \iff h_1(M)=r$. In this case $\pi$ is perfect over $\n F_q[T]$.
\medskip
Anderson t-motives $M$ satisfying these conditions are called uniformizable.
\medskip
There exists a lattice map $$H_1(M)\overset{\al}\to{\to}(\n F_q[\th])^{h_1(M)}\overset{\be}\to{\hookrightarrow}\p^n$$ where $\al$ is an (abstract) isomorphism defined by the condition $\al(T)=\th$ (it serves only in order to identify $T$ and $\th$) and $\be$ is an inclusion of $\n F_q[\th]$-modules, see [G], Section 5.9. The composition inclusion $\be\circ\al$ is an analog of the above $\ga$ for abelian varieties. The image $\be\circ\al(H_1(M))$ is denoted by $L(M)$, it is a $\n F_q[\th]$-lattice of rank $h_1(M)$ in $\p^n$.
\medskip
There is a general
\medskip
{\bf Problem 1.6.} What is a relation between the set of $\n F_q[\th]$-lattices of rank $r$ in $\p^n$, up to $\p$-isomorphisms of $\p^n$, and the set of uniformizable Anderson t-motives of rank $r$ and dimension $n$? Have we some analog of Theorem 1.2?
\medskip
There exists a notion of purity of Anderson t-motives $M$ (see [G], 5.5.2 for the definition). We can expect that in Problem 1.6 we must consider only pure t-motives: conjecturally, the lattice map $M \mapsto L(M)$ has a fibre of dimension $\ge1$ if we consider the set of all Anderson t-motives. Taking into consideration
\medskip
{\bf Theorem 1.7} ([H], Theorem 3.2). The dimension of the moduli space of pure t-motives of rank $r$ and dimension $n$ is equal to $n(r-n)$.
\medskip
and the obvious fact that the moduli space of lattices of rank $r$ in $\p^n$ has the same dimension $n(r-n)$ we can state
\medskip
{\bf Conjecture 1.8.} The image of the lattice map $M \mapsto L(M)$ from the set of pure uniformizable t-motives to the set of lattices is open, and its fibre at a generic point is discrete.
\medskip
The following is known:
\medskip
{\bf Theorem 1.9} (Drinfeld, [Dr]). All t-motives of dimension 1 ( = Drinfeld modules) are pure and uniformizable. There is a 1 -- 1 correspondence between Drinfeld modules of rank $r$ over $\p$ and lattices of rank $r$ in $\p$.
\medskip
For $n=r-1$ the duality theory gives us an immediate corollary of Theorem 1.9:
\medskip
{\bf Corollary 1.10} ([GL07], Corollary 8.4). All pure t-motives of rank $r$ and dimension $r-1$ over $\p$ are uniformizable. There is a 1 -- 1 correspondence between their set, and the set of lattices of rank $r$ in $\p^{r-1}$ having dual.\footnotemark \footnotetext{There is a notion of duality of lattices, see [GL07], Definition 2.3, and [GL07], Section 3. We do not need details here.}
\medskip
Not all such lattices have dual, but almost all, i.e. even in this simple case the correspondence is not strictly 1 -- 1, but only an "almost 1 -- 1".
\medskip
We see that for pure $M$ the minimal values of $r$, $n$ for which we can expect $h^1, \ h_1 <r$ are $r=4$, $n=2$. We shall consider exactly this case. More generally, we shall consider a class of t-motives of dimension $n$ and rank $r=2n$ defined by the equation (2.3) below. We see that they are defined by a matrix $A\in M_{n\times n}(\p)$; the corresponding t-motive is denoted by $M(A)$. All these t-motives are pure, but not all uniformizable.
\medskip
Now we can formulate the results of the present paper. We calculate $h^1$ for all t-motives $M(A)$ where $A$ is of the form (3.4) below, and for some $A$ of the form (3.5) below.

In particular, we describe all uniformizable $M(A)$ where $A$ is of the form (3.4). This is a step to a solution of the Problem 1.6: clearly we need first to describe explicitly the set of all uniformizable Anderson t-motives. Finding (description) of their lattices is a subject of further research.
\medskip
Let us indicate some earlier results to a solution of the Problem 1.6, and some related problems. An explicit description of lattices in $\p^n$ is given in terms of their Siegel matrices $\Cal S$, see for example [GL17], Definition 1.5 (this definition is completely analogous to the definition of Siegel matrices of lattices in $\n C^g$).
\medskip
The main result of [GL17] is, roughly speaking, the 1 -- 1 correspondence between the set of t-motives $M(A)$ defined by (2.3) whose $A$ is in a neighborhood of 0, and the set of lattices whose $\Cal S$ is in a neighborhood of a fixed Siegel matrix $\Cal S_0$ (the main difficulty is to show that we have the same action of some groups on the set of $\Cal S$ and $A$). In particular, it is shown that if all entries $a_{ij}$ of $A$ satisfy $$\ord a_{ij}> \frac q{q^2-1}\eqno{(1.11)}$$ then $M(A)$ is uniformizable ([GL17], end of page 383 and Proposition 2). There is a natural problem to improve this estimate, i.e. to answer
\medskip
{\bf Open problem 1.12.} Let $n, \ q$ be arbitrary. What is the minimal value of $C(n,q)$ satisfying the property: If $A$ from (2.3) satisfies $\forall \ i, \ j$ \ $\ord a_{ij}>C(n,q)$ (version strict inequality) or $\ord a_{ij}\ge C(n,q)$ (version non-strict inequality) then $M(A)$ uniformizable.
\medskip
The above result means $C(n,q)\le \frac q{q^2-1}$. Theorem 1.9 implies $C(1,q)=-\infty$. [G], Example 5.9.9 means $C(n,q) \ge -\frac{q^2}{q-1}$ for $n\ge2$, see Remark 6.3. If we restrict ourselves by $A$ of the form (3.4) then Proposition 6.2 gives $C_{(3.4)}(2,q) =-\frac{q^2}{q-1}$, version strict inequality (here $C_{(3.4)}(2,q)$ means the minimal value of the above $C(n,q)$ for matrices of the form (3.4)).
\medskip
The main result of [GL18] is finding of an explicit method (solution of an affine equation, see (3.1)) for calculation of $h^1(M)$ where $M$ belongs to the same set of $M(A)$, case $n=2$. There exist matrices $A$ such that the application of this method to these $A$ permits to show that not always $h^1(M)=h_1(M)$.
\medskip
Methods of [GL17], [GL18] are essentially different. In [GL17] we find explicitly $L(M)$ which is the kernel of the exponential map of $M$ (see [G], Section 5.9) using a method of successive approximations. This method can be applied to $M$ defined by (2.3), case any $n$, but only for $A$ sufficiently close to 0. For all these $A$ we have $h_1(M(A))=h^1(M(A))=2n$. In [GL18] we calculate $h^1(M)$ solving explicitly an affine equation = a system of polynomial equations, see (3.1). This method can be applied for any $A$, but only for $n=2$ (for $n>2$ the calculations seem to be too difficult).
\medskip
{\bf 2. Definitions on Anderson t-motives.}
\medskip
{\bf Definition 2.1.} The Anderson ring $\p[T,\tau]$ is
the ring of non-commutative polynomials in two variables $T$, $\tau$ over $\p$ satisfying the following relations (here
$a \in \p$):
$$Ta=aT, \ T\tau = \tau T, \ \tau a = a^q \tau $$
Subrings of $\p[T,\tau]$ generated by $\tau$, resp. $T$ are denoted by $\p\{\tau\}$ (a ring of non-commutative polynomials in one variable), resp. $\p[T]$ (the ordinary ring of (commutative) polynomials in one variable).
\medskip
{\bf Definition 2.2.} ([G], 5.4.2, 5.4.18, 5.4.16). A t-motive\footnotemark \footnotetext{Terminology of Anderson; Goss calls these objects abelian t-motives.} $M$ is a left
$\p[T, \tau]$-module which is free and finitely generated as both $\p[T]$-,
$\p\{\tau\}$-module and such that
$$ \exists \goth m= \goth m(M) \ \hbox{such that}\ (T-\theta)^ \goth m M/\tau M=0\eqno{(2.2.1)}$$
\medskip
We shall consider only t-motives for which $\g m=1$.
\medskip
{\bf 2.2.2.} The dimension of $M$ over $\p\{\tau\}$ (resp. $\p[T]$) is denoted by $n$ (resp.
$r$), these numbers are called the dimension and rank of $M$.
\medskip
Let $e_*=(e_1, ..., e_n)^t$ (here and below $t$ means transposition) be the vector column of elements of a basis
of $M$ over $\p\{\tau\}$. To define $M$, it is sufficient to define the multiplication by $T$ of $e_1, ..., e_n$. We shall consider $M$ where the multiplication by $T$ is given by the formula:

$$Te_*=\theta e_*+A\tau e_* +\tau^2e_*\eqno{(2.3)}$$

where $A\in M_{n\times n}(\p)$. The Anderson t-motive defined by (2.3) is denoted by $M(A)$. It is pure of dimension $n$ and rank $2n$.
We shall consider only the case $n=2$. For this case (2.3) becomes

$$T\left(\matrix e_1\\ e_2 \endmatrix \right)=\theta\left(\matrix e_1\\ e_2 \endmatrix \right)+\left(\matrix a_{11}& a_{12} \\  a_{21}& a_{22} \endmatrix \right)\tau \left(\matrix e_1\\ e_2 \endmatrix \right) +\tau^2\left(\matrix e_1\\ e_2 \endmatrix \right)\eqno{(2.4)}$$
where $A=\left(\matrix a_{11}& a_{12} \\  a_{21}& a_{22} \endmatrix \right)\in M_{2\times2}(\p)$ is a matrix.
\medskip
We have $M(A)'=M(A^t)$ (see [GL18], Lemma 1.10.2).
\medskip
Since our purpose is to find $h^1(M)$, we repeat a definition of $H^1(M)$ here. First, we denote by $\p\{T\}$ a subring of $\p[[T]]$ formed by series $\sum_{i=0}^\infty a_iT^i$ such that lim $a_i=0 \ ( \iff \ord a_i\to +\infty$). $\tau$ acts on $\p\{T\}$ by the formula $\tau(\sum_{i=0}^\infty a_iT^i)=\sum_{i=0}^\infty a_i^qT^i$.
\medskip
Now, we define
$$M\{T\}:=M\otimes_{\p[T]}\p\{T\}$$

$\tau$ acts on $M\{T\}$ by the standard formula of the action of an operator on tensor product: $\tau(\al\otimes\be)=\tau(\al)\otimes\tau(\be)$ (see [G], 5.9.11.1).
\medskip
{\bf Definition 2.5.} $H^1(M)=M\{T\}^\tau$ (the set of $\tau$-stable elements).
\medskip
{\bf 2.6.} This definition should be understood as follows. We denote $M[[T]]:=M\otimes_{\p[T]}\p[[T]]$ with the $\tau$-action on the tensor product. We have $M\{T\}^\tau$ is $M[[T]]^\tau \cap M\{T\}$, i.e. $H^1(M)$ is the set of $\tau$-invariant series whose coefficients tend to 0.
\medskip
{\bf 3. Affine equations: definitions and results of [GL18].} For general definitions concerning the affine equations see [GL18], Section 2. Here we repeat them for a particular case that we need. Let $A=\left(\matrix a_{11}& a_{12} \\  a_{21}& a_{22} \endmatrix \right)$ be from (2.4) above. We associate it 9 numbers $a_4, \dots, a_0, b_{14}, b_{13}, b_{12}, b_{24}$ as follows (see [GL18], (3.9), case $\ve=0, \ a_{21}\ne0$):
$$a_4=\frac{\th^{q^3+q^2}}{a_{21}^{q^2}}; \ \ \ a_3=\frac{a_{11}^{q^2}\th^{q^2}}{a_{21}^{q^2}}+\frac{a_{22}^q\th^{q^2}}{a_{21}^q}; \ \ \ a_2=\frac{\th^q}{a_{21}}+\frac{\th^{q^2}}{a_{21}^{q^2}}+\frac{a_{11}^qa_{22}^q}{a_{21}^q}-a_{12}^q; $$ $$ a_1=\frac{a_{11}}{a_{21}}+\frac{a_{22}^q}{a_{21}^q}; \ \ \ a_0=\frac1{a_{21}};\eqno{(3.0)}$$ $$b_{14}=-\frac{\th^{q^3}+\th^{q^2}}{a_{21}^{q^2}}; \ \ \ b_{13}=-\frac{a_{11}^{q^2}}{a_{21}^{q^2}}-\frac{a_{22}^q}{a_{21}^q}; \ \ \ b_{12}=-\frac{1}{a_{21}}-\frac{1}{a_{21}^{q^2}}; \ \ \ b_{24}=\frac{1}{a_{21}^{q^2}}. $$
The affine equation corresponding to $M(A)$ is a series of equations with unknowns $x_0, x_1, x_2, \dots$. The $i$-th equation ($i=0,1,2,\dots)$ of the series has the form (here $x_\al=0$ for $\al<0$)
$$a_4x_i^{q^4}+a_{3}x_i^{q^{3}}+a_2x_i^{q^2}+a_1x_i^q+a_0x_i \ + \ b_{14}x_{i-1}^{q^4}+b_{13}x_{i-1}^{q^3}+b_{12}x_{i-1}^{q^2}+b_{24}x_{i-2}^{q^4}=0\eqno{(3.1)}$$
The terms $a_{k}x_i^{q^{k}}$ ($k=0,\dots,4$) are called the head terms of the equations, other terms are called the tail terms. Explicitly, for $i=0$ the equation is
$$a_4x_0^{q^4}+a_{3}x_0^{q^{3}}+a_2x_0^{q^2}+a_1x_0^q+a_0x_0 =0\eqno{(3.2.0)}$$ (only the head terms);
for $i=1$ $$a_4x_1^{q^4}+a_{3}x_1^{q^{3}}+a_2x_1^{q^2}+a_1x_1^q+a_0x_1 \ + \ b_{14}x_{0}^{q^4}+b_{13}x_{0}^{q^3}+b_{12}x_{0}^{q^2}=0\eqno{(3.2.1)}$$
for $i=2$ $$a_4x_2^{q^4}+a_{3}x_2^{q^{3}}+a_2x_2^{q^2}+a_1x_2^q+a_0x_2 \ + \ b_{14}x_{1}^{q^4}+b_{13}x_{1}^{q^3}+b_{12}x_{1}^{q^2}+b_{24}x_{0}^{q^4}=0 \eqno{(3.2.2)}$$
for $i=3$ $$a_4x_3^{q^4}+a_{3}x_3^{q^{3}}+a_2x_3^{q^2}+a_1x_3^q+a_0x_3 \ + \ b_{14}x_{2}^{q^4}+b_{13}x_{2}^{q^3}+b_{12}x_{2}^{q^2}+b_{24}x_{1}^{q^4}=0 \eqno{(3.2.3)}$$
etc., we have analogous formulas $(3.2.i)$.
\medskip
The set of solutions to (3.2.0) is a 4-dimensional $\n F_q$-vector space in $\p$ denoted by $S_0$. For any such fixed solution $x_0$ the set of $x_1$ satisfying (3.2.1) with the given $x_0$ is an affine space over $S_0$; the same holds for all subsequent equations. This explains the terminology (affine equations).
\medskip
Let $x_0, x_1, x_2, \dots$ be a solution to (3.1). We associate it an element $x_0+x_1T+x_2T^2+... \in\p[[T]]$. The set of these elements is a free $\n F_q[[T]]$-module of rank 4. A solution $x_0, x_1, x_2, \dots$ (and the associated element $x_0+x_1T+x_2T^2+... $) is called a small solution if $\lim  x_i=0$ ($\iff \lim \ord x_i=+\infty$). The set of small solutions is a free $\n F_q[T]$-module.
\medskip
{\bf Theorem 3.3.} $h^1(M(A))$ is the rank of the $\n F_q[T]$-module of small solutions to (3.1).
\medskip
This follows from the calculations of [GL18], Section 3. Non-formally, the meaning of this theorem is the following. Any (non necessarily small) solution to the affine equation (3.1) corresponds to an element of $(M\otimes_{\p[T]}\p[[T]])^\tau$. Condition that a solution is small is equivalent to the condition that it belongs to $H^1(M)=(M\otimes_{\p[T]}\p\{T\})^\tau$.
\medskip
So, our general problem is to solve (3.1) and to find $h^1(A)$ for all $A\in M_{2\times2}(\p)$. This is too complicated, see below, so we consider only some particular cases. Namely, in Section 3 we solve this problem completely for $A$ having the form $$A=\left(\matrix 0& a_{12} \\  a_{21}& 0 \endmatrix \right)\eqno{(3.4)}$$ (because the matrix $A$ from [G], Example 5.9.9 having $h^1=0$ is of this form (after some $\p$-linear change of basis, see Remark 6.3)), and in Section 5 we start to solve this problem for the matrices $A$ having the form $$\left(\matrix a_{11}& a_{12} \\  a_{21}& 0 \endmatrix \right)\hbox{ where $\ord a_{11}=-1$}\eqno{(3.5)}$$ (because the matrix $A$ from [GL18], Section 4 --- giving a counterexample to $h^1=h_1$ --- is of this form).
\medskip
{\bf 3.6.} Our purpose is to answer some questions of [GL18], Section 0.3. As the first step, we want to find all possible values of $h^1(A)$. There are examples $h^1=0,1,4$, but until now there is no examples $h^1=2$ or 3. Particularly,
it turns out that $h^1$ of all $M(A)$ of Section 3 is equal to 0 or 4.
\medskip
Hence, at the moment we are far from a solution of questions of [GL18], Section 0.3, for example: what are possible triples of numbers \{$h^1(M), \ h_1(M)$, rank of pairing $\pi$\} ? Surprisingly, we have no other non-trivial pairings except the ones that are covered by Theorem 1.5. Moreover,
it seems that for "almost all" matrices $A$ we have $h^1(M(A))=4$. The words "almost all" should be understood as follows. Let us consider the ord$_4$ map from $M_{2\times2}(\p)$ to $(\n Q\cup+\infty)^4$ defined as ord's of entries of $A$: $$\hbox{ord}_4(A)= (\ord a_{11}, \ \ord a_{12}, \ \ord a_{21}, \ \ord a_{22})\in (\n Q\cup+\infty)^4$$
\medskip
{\bf Conjecture 3.7.} There exists a subset $\g U$ of $(\n Q\cup+\infty)^4$ which is the complement to a union of countably many (maybe even finitely many - we do not know) linear subspaces of dimension $\le3$ such that if ord$_4(A)\in \g U$ then $h^1(M(A))=4$.
\medskip
For $A$ of the form (3.4) this conjecture is confirmed by Proposition 6.1. Moreover, results of some computer calculations (finding of minimal chains, see 4.4 below) made by the authors also support this conjecture.
\medskip
{\bf 4. Method of calculation.} We shall solve consecutively for $i=0,1,2,\dots$ the equations (3.1) for a given $A$. A solution $(x_0, \ x_1, \ x_2,\dots)$ will be denoted by $\{x\}$; if $\{x\}$ carries a subscript $\{x_k\}$ then the corresponding $(x_0, \ x_1, \ x_2,\dots)$ will be denoted as $(x_{k0}, \ x_{k1}, \ x_{k2},\dots)$. There is a trivial
\medskip
{\bf Lemma 4.1.} ( = [GL18], Proposition 2.3.) Let $\{x_1\},\dots, \{x_4\}$ be solutions to (3.1). This set is a basis of the set of solutions to (3.1) over $\n F_q[[T]]$ iff $x_{10},\dots,x_{40}$ is a basis of $S_0$ over $\n F_q$. $\square$
\medskip
To calculate $\ord x_i$ we use the notion of Newton polygon. Let us give the corresponding definitions. Let $P=\sum_{i=0}^n c_ix^i$ be a polynomial, $c_i\in\p$. We associate it a set of $n+1$ points on a plane whose coordinates are $(i, \ord c_i)$, $i=0,\dots,n$. These points are called the Newton points of $P$. The Newton polygon of $P$ is the lower convex hull of its Newton points. The ord's of the roots of $P$ are the minus slopes of the segments of its Newton polygon.
\medskip
Let us apply this notion to an affine equation. The $0$-th polynomial of an affine equation gives us a set of Newton points $(1, \ord a_0), \ (q,\ord a_1), ...,\ (q^4, \ord a_4)$. If the Newton polygon of the $0$-th polynomial consists of 4 segments then their minus slopes are the ord's of the four elements of a basis of $S_0$ over $\n F_q$. According the above notations, we denote these basis elements by $x_{10},\dots,x_{40}$. If the Newton polygon of the $0$-th polynomial consists of less than 4 segments then some of ord's of $x_{10},\dots,x_{40}$ are equal.
\medskip
Now let us consider the $i$-th affine equation. Its unknown is $x_i$, we consider $x_0,\dots,x_{i-1}$ as already known (fixed). We denote the ord of the sum of the tail terms of the $i$-th equation by $g_i$, hence the Newton points of the $i$-th equation are $(0,g_i), \ (1, \ord a_0), \ (q,\ord a_1), ...,\ (q^4, \ord a_4)$.
\medskip
Let us recall a notion of a minimal solution from [GL18], Definition 2.5: a solution to (3.1) $\{x\}=(x_0, \ x_1, \ x_2,\dots)$ is called a minimal solution (generated by $x_0$) if $\forall \ i>0$ it satisfies the following condition: $\ord x_i$ corresponds to the leftmost segment of the Newton polygon of $(3.2.i)$ ( = the $i$-th equation of (3.1)), i.e. $\ord x_{i}$ has the maximal possible value amongst ord's of solutions to the $i$-th affine equation for fixed $x_0, x_1, x_2, \dots, x_{i-1}$.
\medskip
Clearly a minimal solution generated by $x_0$ is not unique. Moreover, it can happen that even $\ord x_i$ of a minimal solution vary ($i$ is fixed, minimal solutions vary): this can occur if there exists a jump of valuation\footnotemark \footnotetext{A jump of valuation is a situation when $\ord a = \ord b$ and $\ord a+b> \ord a$.} of the tail terms.
\medskip
Conversely, let us formalize the situation when ord's of tail terms are different. We repeat [GL18], Definition 2.4: a solution $\{x\}=(x_0, \ x_1, \ x_2,\dots)$ is called simple  if for all $i$ we have: among all tail members of the $i$-th equation $(3.2.i)$ for this $\{x\}$ (i.e. obtained while we substitute $x_0,x_1,...,x_{i_0-1}$) there exists only one term whose ord is the minimal one.
\medskip
Clearly if one minimal solution generated by $x_0$ is simple then all minimal solutions generated by $x_0$ are also simple, and for any $i$ the $\ord x_i$ are the same for all minimal solutions generated by $x_0$. They depend only on $i$ and on ord's of $a_0,\dots,a_4, \ b_{12}, \ b_{13}, \ b_{14}, \ b_{24}$.
\medskip
{\bf 4.2.} Further, we can define a sequence of ord's (called a simple minimal sequence) as follows. We start from $\ord x_{i0}$ for some fixed $i$, $i=1,\dots,4$. We define $\tilde g_1$ as the minimum of the ord's of the tail terms of (3.2.1). We define $v_1$ (abbreviation of "valuation of $x_{i1}$") as the minus slope of the leftmost segment of the convex hull of $(0,\tilde g_1), \ (1, \ord a_0), \ (q,\ord a_1), ...,\ (q^4, \ord a_4)$. Now we consider the ord's of the tail terms of (3.2.2) where instead of $\ord x_1$ we substitute $v_1$. Exactly, we define $\tilde g_2:= \min(\ord b_{1k} + q^k v_1 $, $k=2,3,4$, and $\ord b_{24} + q^4 \ord x_{i0})$. We define $v_2$ as the minus slope of the leftmost segment of the convex hull of $(0,\tilde g_2), \ (1, \ord a_0), \ (q,\ord a_1), ...,\ (q^4, \ord a_4)$. Continuing this process we define $\tilde g_3:= \min(\ord b_{1k} + q^k v_2 $, $k=2,3,4$, and $\ord b_{24} + q^4 v_1)$, we define $v_3$ etc.
\medskip
For any minimal solution $x_{i0}, x_{i1}, x_{i2}, \dots$ generated by $x_{i0}$ the following inequality holds: $$\forall \ j\ \ \ \ord x_{ij}\ge v_j$$ (it is proved immediately by induction by $j$). In particular if $v_j\to +\infty$ then a minimal solution generated by $x_{i0}$ is a small solution. Moreover we have
\medskip
{\bf Proposition 4.3.} If
$\ord x_{10}\ge...\ge \ord x_{40}$ and the sequence $v_j$ for $x_{40}$ tends to $+\infty$ then $h^1=4$. $\square$
\medskip
{\bf 4.4.} This simple criterion was used in computer calculations. We chose random values of $\ord a_{ij}$, $i, j=1,\ 2$ (i.e. 4 random rational numbers). We assumed that there were no jumps of valuation in terms of the formulas (3.0), i.e. we put ord of $a_0,\dots, a_4,\ b_{12},\ b_{13},\ b_{14},\ b_{24}$ as the minimum of the ord's of their terms. Further, we calculated (according the algorithm described in (4.2) ) the values of $v_j$ starting from $\ord x_{40}$ (the minimal value, see 4.3). In all considered cases we had $v_j\to+\infty$. This supports Conjecture 3.7.
\medskip
{\bf 4.5.} It is easy to see that if the structure of the convex hull of Newton points does not depend on $i$ then the behavior of $v_i$ is either linear or exponential:
\medskip
$v_i= \al + \be i$ (linear, $\al, \ \be$ are constants), or
\medskip
$v_i=\al + \be\cdot\ga^i$ (exponential, $\al, \ \be, \ \ga>0$ are constants). It can happen $\be>0$ or $\be<0$, $\ga>1$ or $\ga<1$.
\medskip
{\bf 4.6.} Proofs of all propositions of the present paper are similar. We use induction to find the above $\al, \ \be, \ \ga$ and according the cases $\ga>1$ or $\ga<1$ we get the result. Difficulties occur if there are different types of complex hulls, or jumps of valuation.
\medskip
{\bf 4.7. Further research.} The reader will see that in order to solve the problem of calculation of $h^1$ completely (for all $A$) we must use computers. There are questions:
\medskip
{\bf 4.7.1.} Is it true that for any fixed $A$ and for all minimal solutions there exists $i_0$ such that for $i>i_0$ the types of convex hulls of the Newton polygons of $(3.2.i)$ are the same (and hence $v_i$ behave like in (4.5)) ?
\medskip
{\bf 4.7.2.} Does exist an algorithm of calculation of $h^1(M(A))$ for any fixed $A$?
\medskip
Even if the answer to (4.7.1) is YES, we cannot guarantee existence of such algorithm, because a linear combination of non-small solutions can give a small solution, i.e. $h^1(M(A))$ can be bigger than the quantity of independent small minimal solutions.
\medskip
{\bf 4.7.3.} Is the quantity of types of the matrices $A$ (see below for the types) finite or not?
\medskip
In order to prove or to disprove Conjecture 1.8, we should study the lattices corresponding to uniformizable t-motives. The first question in this direction is
\medskip
{\bf Problem 4.7.4.} Find (describe) Siegel matrices corresponding to $M(A)$ where $A$ is of the form (3.4) such that $M(A)$ is uniformizable.
\medskip
{\bf 5. Case of $A$ of the form (3.4).} We use notations $v$, $t$ where $v=\ord a_{21}$ and $t:=\ord a_2=\ord(\frac{\th^q}{a_{21}}+\frac{\th^{q^2}}{a_{21}^{q^2}}-a_{12}^q)$. Since $a_{12}$ appears in the formulas (3.0) only once (in $a_2$) we get that $a_2$ defines $a_{12}$ uniquely.
\medskip
{\bf Theorem 5.1.} For $t\ge-v, \ v\le-\frac{q^2}{q-1}$ and for $t\ge-q^2v, \ v\ge\frac{1}{q-1}$ we have $h^1(M(A))=0$, for other values of $t, \ v$ we have $h^1(M(A))=4$.
\medskip
{\bf Remark.} The set of points on $(v,t)$-plane having $h^1=0$ is the union of interiors and borders of two angles (denoted by $\Cal A_-$, resp. $\Cal A_+$ according the sign of $v$) whose intersection is empty.
\medskip
{\bf Proof.} It will consist of the below Propositions 5.7, 5.9, 5.12. Formulas (3.0) for the present case give us
\medskip
$\ord a_4=-q^2v-q^3-q^2$; \ \ $a_3=0$; \ \ $\ord a_2=t$; \ \ $a_1=0$; \ \ $\ord a_0=-v$;
\medskip
$\ord b_{14}=-q^2v-q^3$; \ \ $b_{13}=0$; \ \ $\ord b_{12}=-\de v$ where $\de=1$ for $v<0$, $\de=q^2$ for $v>0$, and $\ord b_{12}\ge0$ for $v=0$; $\ord b_{24}=-q^2v$.
\medskip
We denote vertices of the angles $\Cal A_-$, resp. $\Cal A_+$ by $V_-$, resp. $V_+$.
We have:
$$t<-\frac{2q^2}{q^2+1}v-\frac{q^3+q^2}{q^2+1}\eqno{(5.2)}$$ is the condition that the Newton polygon of (3.2.0) consists of two segments, their vertices are $(1,-v)$, $(q^2, t)$, $(q^4, -q^2v-q^3-q^2)$.
\medskip
{\bf 5.3.} We have: the straight line $t=-\frac{2q^2}{q^2+1}v-\frac{q^3+q^2}{q^2+1}$ passes through $V_-$, $V_+$.
\medskip
If (5.2) holds then the ord's of $x_{j0}$ are:
$$\ord x_{10}=\ord x_{20}=-\frac1{q^2-1}t-\frac1{q^2-1}v\eqno{(5.4)}$$ $$\ord x_{30}=\ord x_{40}= \frac1{q^4-q^2}t+\frac1{q^2-1}v+\frac1{q-1}\eqno{(5.5)}$$
\medskip
{\bf Case 5.6.} $t\le -\frac{2q^2}{q^2+1}v-\frac{q^3+q^2}{q^2+1}$ and $t\ge-\de v$. (5.3) implies that the set of $(v,t)$ satisfying these conditions is the union of two angles (denoted by $\Cal A_{1-}, \ \Cal A_{1+}$) such that $\Cal A_{1-}\subset \Cal A_-$, $\Cal A_{1+}\subset \Cal A_+$, their vertices are $V_-$, resp. $V_+$, and they have common rays with $\Cal A_-$, resp. $\Cal A_+$. Let us prove that in this case $h^1=0$. We use notations of Section 4: $x_0, \ x_1, \ x_2, \dots$ is any solution to the affine equation (3.1), and for any fixed $i>0$ let $g_i:=\ord(b_{14}x_{i-1}^{q^4}+b_{13}x_{i-1}^{q^3}+b_{12}x_{i-1}^{q^2}+b_{24}x_{i-2}^{q^4})$.
\medskip
{\bf Proposition 5.7.} For the case 5.6 for any solution $x_0, \ x_1, \ x_2, \dots$ we have:
\medskip
(5.7.1) $\forall \ i$ \ \ $\frac1{q^2}v+\frac1{q-1}\le \ord x_i\le -\frac1{q^2-1}t-\frac1{q^2-1}v$ if $v<0$,
\medskip
(5.7.2) $\frac1{q-1}\le \ord x_i\le -\frac1{q^2-1}t-\frac1{q^2-1}v$ if $v>0$;
\medskip
(5.7.3) $\forall \ i$ \ \ $\ord x_i\le \ord x_{i-1}$.
\medskip
{\bf Proof.} First, we check that conditions of Case 5.6 imply
\medskip
(5.7.4) $\frac1{q^2}v+\frac1{q-1}\le -\frac1{q^2-1}t-\frac1{q^2-1}v$ if $v<0$,
\medskip
(5.7.5) $\frac1{q-1}\le -\frac1{q^2-1}t-\frac1{q^2-1}v$ if $v>0$.
\medskip
Case $v<0$: (5.7.4) is $\frac{-2q^2+1}{q^2}v-q-1\ge t$. Because of (5.2), it is sufficient to prove that $\frac{-2q^2+1}{q^2}v-q-1\ge -\frac{2q^2}{q^2+1}v-\frac{q^3+q^2}{q^2+1}$. This is equivalent $v\le -\frac{q^2}{q-1}$, i.e. (5.7.4) holds. Case $v>0$: (5.7.5) is $t\le -v-q-1$, this holds for Case 5.6, $v>0$.
\medskip
Now we use induction by $i$. Let $i=0$. Condition $t\le -\frac{2q^2}{q^2+1}v-\frac{q^3+q^2}{q^2+1}$ implies that $\ord x_0= -\frac1{q^2-1}t-\frac1{q^2-1}v$ or $\ord x_0=\frac1{q^4-q^2}t+\frac1{q^2-1}v+\frac1{q-1}$, see (5.4), (5.5). Both these values satisfy (5.7.1-2), if the conditions of the Case 5.6 hold.
\medskip
Now we assume that (5.7.1-3) hold for a fixed value of $i$ (where for $i=0$ condition (5.7.3) is omitted), and prove that they hold for $i+1$. First, let us prove that (5.7.1-3) imply that $$\forall \ i \ \ \ \ord b_{14}x_{i}^{q^4}>\ord b_{12}x_{i}^{q^2}<\ord b_{24}x_{i-1}^{q^4}$$ This is immediate. For $v<0$ we have $\ord b_{14}x_{i}^{q^4}>\ord b_{12}x_{i}^{q^2}$ is equivalent to $\ord x_i>\frac1{q^2}v+\frac{q}{q^2-1}$, which holds because of (5.7.1).
For $v>0$ we have $\ord b_{14}x_{i}^{q^4}>\ord b_{12}x_{i}^{q^2}$ is equivalent to $\ord x_i>\frac{q}{q^2-1}$, which holds because of (5.7.2).
\medskip
Further, we have $$\ord b_{24}x_{i-1}^{q^4}=-q^2v+q^4\ord x_{i-1}\ge -q^2v+q^4\ord x_{i}>-\de v+q^2\ord x_{i}=\ord b_{12}x_{i}^{q^2}$$ (the inner inequality because of $\ord x_{i}>\frac1{q^2}v$ for $v<0$, $\ord x_{i}>0$ for $v>0$).
\medskip
 Hence, we have $g_{i+1}=\ord b_{12}x_{i}^{q^2}=-\de v+q^2\ord x_{i}$. (5.7.1), (5.7.2) imply $$\frac{q^2}{q-1}\le g_{i+1}\le -\frac{q^2}{q^2-1}t-\frac{2q^2-1}{q^2-1}v \hbox{ (case $v<0$),}$$ $$q^2-q^2v\le g_{i+1}\le -\frac{q^2}{q^2-1}t-\frac{q^4}{q^2-1}v \hbox{ (case $v>0$)}$$
The line defined by the points $(q^2, t)$, $(q^4, -q^2v-q^3-q^2)$ (the second segment of the Newton polygon) crosses the $t$-axis at the point $$(0, \frac{q^2}{q^2-1}t+\frac{q^2}{q^2-1}v+\frac{q^2}{q-1})\eqno{(5.8)}$$ and the line defined by points $(1, -v)$, $(q^2, t)$ (the first segment of the Newton polygon) crosses the $t$-axis at the point $$(0, -\frac1{q^2-1}t-\frac{q^2}{q^2-1}v)$$ The fact that for $(t,v)$ defined by the conditions of the Case 5.6 we have the inequalities $$-\frac{1}{q^2-1}t-\frac{q^2}{q^2-1}v\ge -\frac{q^2}{q^2-1}t-\frac{2q^2-1}{q^2-1}v \hbox{ (case $v<0$)}$$
$$-\frac{1}{q^2-1}t-\frac{q^2}{q^2-1}v\ge -\frac{q^2}{q^2-1}t-\frac{q^4}{q^2-1}v \hbox{ (case $v>0$)}$$
implies that the convex hull of the point $(0,g_{i+1})$ and the Newton polygon of the head of the equation (3.1) consists of the points $$\hbox{ $(0,g_{i+1})$, $(q^2, t)$, $(q^4, -q^2v-q^3-q^2)$ if $g_{i+1}>\frac{q^2}{q^2-1}t+\frac{q^2}{q^2-1}v+\frac{q^2}{q-1}$}$$ $$\hbox{ $(0,g_{i+1})$, $(q^4, -q^2v-q^3-q^2)$ if  $g_{i+1}\le\frac{q^2}{q^2-1}t+\frac{q^2}{q^2-1}v+\frac{q^2}{q-1}$}$$ Hence, if $g_{i+1}> \frac{q^2}{q^2-1}t+\frac{q^2}{q^2-1}v+\frac{q^2}{q-1}$ then there are two values of $\ord x_{i+1}$: $$\ord x_{i+1,1}=\ord x_{i} -\frac1{q^2}t-\frac\de{q^2}v\hbox{ and }\ord x_{i+1,2}= \frac1{q^4-q^2}t+\frac1{q^2-1}v+\frac1{q-1}$$ where $\ord x_{i+1,1}>\ord x_{i+1,2}$.
\medskip

If $g_{i+1}\le \frac{q^2}{q^2-1}t+\frac{q^2}{q^2-1}v+\frac{q^2}{q-1}$ then there is one value of $\ord x_{i+1}$: $$\ord x_{i+1}=\frac1{q^2}\ord x_{i}+\frac{q^2-\de}{q^4}v+\frac{q+1}{q^2}$$ In all cases we have (5.7.1-3) for $i$ imply (5.7.1-3) for $i+1$: this is checked immediately. $\square$
\medskip
Proposition 5.7 implies that for the domain $t\le -\frac{2q^2}{q^2+1}v-\frac{q^3+q^2}{q^2+1}$ and $t\ge-\de v$ we have $h^1=0$.
\medskip
{\bf Proposition 5.9}. If $t<-\frac{2q^2}{q^2+1}v-\frac{q^3+q^2}{q^2+1}$, $t<-\de v$ then $h^1=4$.
\medskip
{\bf Proof.} We consider a minimal chain generated by $x_{40}$, see 4.3. We abbreviate $\g A:=\frac{q^2}{q^2-1}t+\frac{q^2}{q^2-1}v+\frac{q^2}{q-1}$ ($\g A$ from (5.8)) and $\al:=\min(\ord b_{14} x_{40}^{q^4}, \ \ord b_{12} x_{40}^{q^2})-\g A$ where $\ord x_{40}$ is from (5.5). Condition $t<-\de v$ implies $\al>0$. We shall prove by induction that for all $i$
$$g_i\ge \g A+\al i\eqno{(5.9.1)}$$
$$\ord x_{4i}\ge \ord x_{40}+\frac\al{q^2}i\eqno{(5.9.2)}$$

(5.9.2) holds for $i=0$ and (5.9.1) holds for $i=1$. Further, truth of (5.9.1) for $i$ implies truth of (5.9.2) for $i$ (if the convex hull of $(0,g_i); \ (1,-v); \ (q^2,t)$ is the segment $(0,g_i); \ (q^2,t)$ and $g_i= \g A+\al i$ then in (5.9.2) for $i$ we have equality; if not then $\ord x_{4i}$ is higher).

So, we should prove that truth of (5.9.2) for $i$ implies truth of (5.9.1) for $i+1$. We have: truth of (5.9.1) for $i+1$ is equivalent to

$$\ord b_{14} x_{4i}^{q^4}\ge \g A+\al(i+1)$$
$$\ord b_{12} x_{4i}^{q^2}\ge \g A+\al(i+1)$$
$$\ord b_{24} x_{4,i-1}^{q^4}\ge \g A+\al(i+1)$$

All these inequalities follow immediately from (5.9.2) for $i$ and the conditions of the proposition.

Finally, (5.9.2) shows that $x_{j0}$ generates a small solution for $j=3,4$. For $j=1,2$ the same inequalities (5.9.1), (5.9.2) also hold, hence the proposition (also, we can use 4.3). $\square$
\medskip
{\bf Case 5.10.} $t\ge-\frac{2q^2}{q^2+1}v-\frac{q^3+q^2}{q^2+1}$. This is a condition that the Newton polygon of the head of the equation consists of one segment, its vertices are $(1,-v)$, $(q^4, -q^2v-q^3-q^2)$, and
\medskip
$$\forall \ j\ \ \ \ord x_{j0}= \frac1{q^2+1}v+\frac{q^2}{(q-1)(q^2+1)}\eqno{(5.11)}$$
\medskip
{\bf Proposition 5.12.} For the Case 5.10 we have $h^1=0$ if $v\ge\frac1{q-1}$, $v\le-\frac{q^2}{q-1}$, and $h^1=4$ otherwise.
\medskip
{\bf Proof.} Induction. We consider five cases; the proofs for all these cases are similar. The below equalities (cases (a) - (d) ) and inequalities (case (e) ) hold for all solutions of (3.1), not necessarily for minimal solutions.
\medskip
(a) Case $v\ge\frac1{q-1}$: We have
$$\forall \ j, \ i\ \ \ \ \ \ord x_{ji}= \frac1{(q^2+1)q^{2i}}v+\frac1{q-1}-\frac1{(q^2+1)(q-1)q^{2i}}\eqno{(5.12.a.1)}$$ Really, for $i=0$ this is (5.11); this argument will be also valid for the below cases (b) - (e). The induction step from $i$ to $i+1$: we have
$$\ord b_{14}x_{ji}^{q^4}=(\frac{q^4}{(q^2+1)q^{2i}}-q^2)v+\frac{q^3}{q-1}-\frac{q^4}{(q^2+1)(q-1)q^{2i}}\eqno{(5.12.a.2)}$$
$$\ord b_{12}x_{ji}^{q^2}=(\frac{q^2}{(q^2+1)q^{2i}}-q^2)v+\frac{q^2}{q-1}-\frac{q^2}{(q^2+1)(q-1)q^{2i}}\eqno{(5.12.a.3)}$$
$$\ord b_{24}x_{j,i-1}^{q^4}=(\frac{q^6}{(q^2+1)q^{2i}}-q^2)v+\frac{q^4}{q-1}-\frac{q^6}{(q^2+1)(q-1)q^{2i}}\eqno{(5.12.a.4)}$$
{\bf (5.12.1)} The below arguments hold for all cases (a) --- (d), and they are slightly modified for the case (e):
\medskip
For $v$ of the present case we have
$$\ord b_{12}x_{ji}^{q^2}<\ord b_{14}x_{ji}^{q^4}, \ \ \ord b_{12}x_{ji}^{q^2}<\ord b_{24}x_{j,i-1}^{q^4}\hbox{ (here $i\ge1$)}\eqno{(5.12.2)}$$ We get that for $v$ of the present case the Newton polygon of the equation for $x_{j,i+1}$ is a segment $(0, \ord b_{12}x_{ji}^{q^2}), (q^4, -q^2v-q^3-q^2)$, hence, we get the induction proposition.
\medskip
(b) Case $v\le-\frac{q^2}{q-1}$. Analogs of (5.12.a.1) --- (5.12.a.4) are:
$$\forall \ j, \ i\ \ \ \ \ord x_{ji}=(\frac1{q^2}-\frac1{(q^4+q^2)q^{2i}})v+\frac1{q-1}- \frac1{(q^2+1)(q-1)q^{2i}}\eqno{(5.12.b.1)}$$
$$\ord b_{14}x_{ji}^{q^4}=-\frac{q^2}{(q^2+1)q^{2i}}v+\frac{q^3}{q-1}-\frac{q^4}{(q^2+1)(q-1)q^{2i}}\eqno{(5.12.b.2)}$$
$$\ord b_{12}x_{ji}^{q^2}=-\frac{1}{(q^2+1)q^{2i}}v+\frac{q^2}{q-1}-\frac{q^2}{(q^2+1)(q-1)q^{2i}}\eqno{(5.12.b.3)}$$
$$\ord b_{24}x_{j,i-1}^{q^4}=-\frac{q^4}{(q^2+1)q^{2i}}v+\frac{q^4}{q-1}-\frac{q^6}{(q^2+1)(q-1)q^{2i}}\eqno{(5.12.b.4)}$$
Arguments (5.12.1) are the same.
\medskip
(c) Case $0< v<\frac1{q-1}$. Analogs of (5.12.a.1) --- (5.12.a.4) are:
$$\forall \ j, \ i\ \ \ \ \ord x_{ji}=(1-\frac{q^2}{q^2+1}q^{2i})v+\frac{q^2}{(q^2+1)(q-1)}q^{2i}\eqno{(5.12.c.1)}$$
$$\ord b_{14}x_{ji}^{q^4}=(q^4-q^2-\frac{q^6}{q^2+1}q^{2i})v+\frac{q^6}{(q^2+1)(q-1)}q^{2i}-q^3\eqno{(5.12.c.2)}$$
$$\ord b_{12}x_{ji}^{q^2}=-\frac{q^4}{q^2+1}q^{2i}v+\frac{q^4}{(q^2+1)(q-1)}q^{2i}\eqno{(5.12.c.3)}$$
$$\ord b_{24}x_{j,i-1}^{q^4}=(q^4-q^2-\frac{q^4}{q^2+1}q^{2i})v+\frac{q^4}{(q^2+1)(q-1)}q^{2i}\eqno{(5.12.c.4)}$$
Arguments (5.12.1) are the same.
\medskip
(d) Case $-\frac{q^2}{q-1}< v<0$. Analogs of (5.12.a.1) --- (5.12.a.4) are:
$$\forall \ j, \ i\ \ \ \ \ord x_{ji}=\frac1{q^2+1}q^{2i} v+\frac{q^2}{(q^2+1)(q-1)}q^{2i} \eqno{(5.12.d.1)}$$
$$\ord b_{14}x_{ji}^{q^4}=(\frac{q^4}{q^2+1}q^{2i}-q^2)v+\frac{q^6}{(q^2+1)(q-1)}q^{2i}-q^3\eqno{(5.12.d.2)}$$
$$\ord b_{12}x_{ji}^{q^2}=(\frac{q^2}{q^2+1}q^{2i}-1)v+\frac{q^4}{(q^2+1)(q-1)}q^{2i}\eqno{(5.12.d.3)}$$
$$\ord b_{24}x_{j,i-1}^{q^4}=(\frac{q^2}{q^2+1}q^{2i}-q^2)v+\frac{q^4}{(q^2+1)(q-1)}q^{2i}\eqno{(5.12.d.4)}$$
Arguments (5.12.1) are the same; in order to prove (5.12.2) for this case it is sufficient to check these inequalities for the values $v=0$ and $v=-\frac{q^2}{q-1}$; truth of (5.12.2) for intermediate values of $v$ holds by linearity.
\medskip
(e) Case $v=0$. It is slightly different from the above cases because we have $\ord b_{12}\ge0$ and equalities (5.12.*.1) --- (5.12.*.4) are replaced by inequalities:
$$\forall \ j, \ i\ \ \ \ \ord x_{ji}\ge\frac{q^2}{(q^2+1)(q-1)}q^{2i}\eqno{(5.12.e.1)}$$
$$\ord b_{14}x_{ji}^{q^4}\ge\frac{q^6}{(q^2+1)(q-1)}q^{2i}-q^3\eqno{(5.12.e.2)}$$
$$\ord b_{12}x_{ji}^{q^2}\ge\frac{q^4}{(q^2+1)(q-1)}q^{2i}\eqno{(5.12.e.3)}$$
$$\ord b_{24}x_{j,i-1}^{q^4}\ge\frac{q^4}{(q^2+1)(q-1)}q^{2i}\eqno{(5.12.e.4)}$$
Formulas  (5.12.2) are modified as follows:
$$\ord b_{14}x_{ji}^{q^4}\ge\frac{q^4}{(q^2+1)(q-1)}q^{2i}; \ \ \ord b_{24}x_{j,i-1}^{q^4}\ge \frac{q^4}{(q^2+1)(q-1)}q^{2i}$$
All other arguments are the same.
\medskip
We get that for the cases (c), (d), (e) we have $\ord x_{ji}\to+\infty$ as $i\to \infty$, while for the cases (a), (b) $\ord x_{ji}$ is bounded. This proves the proposition. $\square$
\medskip
Propositions 5.7, 5.9, 5.12 cover the whole plane $(v,t)$, hence we get the Theorem 5.1. $\square$
\medskip
{\bf 6. Description of the sets $h^1=4$, $h^1=0$ in terms of $a_{12}$, $a_{21}$.}
\medskip
Let $u=\ord a_{12}$ and $v=\ord a_{21}$ as above. Here we find the sets of $h^1=4$, $h^1=0$ on $(u,v)$-plane. This will give us an improved estimate (1.11) for the matrices $A=\left(\matrix 0& a_{12} \\  a_{21}& 0 \endmatrix \right)$, see Proposition 6.2.
\medskip
Further, we verify that the symmetry with respect to $a_{12} \longleftrightarrow a_{21}$ really takes place. Although most likely for a generic $A$ the t-motives $M(A)$, $M(A)'$ are not isomorphic\footnotemark \footnotetext{It is easy to prove that a constant matrix of a change of basis (i.e. a matrix whose coefficients belong to $\p$) does not give an isomorphism between $M(A)$ and $M(A)'$. Most likely the same is true for a matrix of a change of basis whose coefficients belong to $\p\{\tau\}$.}, Theorems 1.5 and 5.1 show that always $h^1(M(A))=h^1(M(A'))$. This is not seen beforehand, because formulas (3.0) are not symmetric with respect to the transposition of $A$.
\medskip
{\bf 6.0.} We have $t=\ord a_2=\ord(\frac{\th^q}{a_{21}}+\frac{\th^{q^2}}{a_{21}^{q^2}}-a_{12}^q)$. There are 3 domains on the $(u,v)$-coordinate plane, according the minimality of ord's of the 3 summands of $a_2$. Let us describe them. There are 3 rays (half-lines) $R_1$, $R_2$, $R_3$ having the same initial point $(u,v)=(-\frac {q}{q+1},-\frac {q}{q+1})$.
\medskip
The ray $R_1$ has the equation $v=-\frac {q}{q+1}$, $u\ge -\frac {q}{q+1}$;
\medskip
The ray $R_2$ has the equation $v=-qu-q$, $u\ge -\frac {q}{q+1}$;
\medskip
The ray $R_3$ has the equation $v=-\frac {1}{q}u-1$, $u\le -\frac {q}{q+1}$;
\medskip
The open domain between the rays $R_i$ and $R_j$ is denoted by $D_{ij}$.
\medskip
In Domain $D_{12}$ we have $\ord\frac{\th^q}{a_{21}}<\ord \frac{\th^{q^2}}{a_{21}^{q^2}}, \ \ \ord\frac{\th^q}{a_{21}}<\ord a_{12}^q$, i.e. $t=-v-q$. The ray $t=-v-q$, $v<-\frac {q}{q+1}$ on the $(v,t)$-coordinate plane is entirely in the domain $h^1=4$ (see Theorem 5.1).
\medskip
In Domain $D_{13}$ we have $\ord \frac{\th^{q^2}}{a_{21}^{q^2}} <\ord\frac{\th^q}{a_{21}}, \ \ \ord \frac{\th^{q^2}}{a_{21}^{q^2}} <\ord a_{12}^q$, i.e. $t=-q^2v-q^2$. The ray $t=-q^2v-q^2$, $v>-\frac {q}{q+1}$ on the $(v,t)$-coordinate plane is entirely in the domain $h^1=4$ (see Theorem 5.1).
\medskip
In Domain $D_{23}$ we have $\ord a_{12}^q < \ord \frac{\th^{q^2}}{a_{21}^{q^2}}, \ \ \ord a_{12}^q < \ord\frac{\th^q}{a_{21}}$, i.e. $t=qu$. The image of $D_{23}$ under the map $(u,v) \to (t,v)$ defined by $(u,v) \mapsto (qu,v)$, is the open interior of the angle formed by the rays $t=-v-q$, $v<-\frac {q}{q+1}$ and $t=-q^2v-q^2$, $v>-\frac {q}{q+1}$. It is entirely in the domain $h^1=4$ (see Theorem 5.1).
\medskip
In particular, we get that Conjecture 3.7 holds for the case under consideration:
\medskip
{\bf Proposition 6.1.} The set of points on $(u,v)$-plane such that $M(A)$ can be non-uniformizable is of dimension less then 2 (i.e. "almost all" $M(A)$, $A$ of the form (3.4), are uniformizable). $\square$
\medskip
For the points of the ray $R_1$ we have: $v=-\frac {q}{q+1}$, $t$ can vary. According Theorem 5.1, all these points have $h^1=4$. Moreover, Theorem 5.1 shows that if $-\frac{q^2}{q-1}<v<\frac{1}{q-1}$ then $h^1=4$. This means that the only points on $(u,v)$-coordinate plane where it can happen $h^1=0$ are subrays $\bar R_2$, resp. $\bar R_3$ of $R_2$, resp. $R_3$ having initial points
$$u=-\frac{q^2}{q-1}, \ v=\frac{1}{q-1}\hbox{ for }\bar R_2$$
$$u=\frac{1}{q-1}, \ v=-\frac{q^2}{q-1}\hbox{ for }\bar R_3$$

{\bf Remark.} (1) We see that $\bar R_2$, $\bar R_3$ are symmetric with respect to the symmetry $u \longleftrightarrow v$ as it must be.
\medskip
(2) Clearly not for all $a_{12}, \ a_{21}$ such that their ord's belong to $\bar R_2$, $\bar R_3$ we have $h^1=0$.
\medskip
Further, we get
\medskip
{\bf Proposition 6.2.} Let $A=\left(\matrix 0& a_{12} \\  a_{21}& 0 \endmatrix \right)$. If $\ord a_{12}, \ \ord a_{21} > -\frac{q^2}{q-1}$ then $M(A)$ is uniformizable. $\square$
\medskip
{\bf  Remark 6.3.} Example [G], 5.9.9 shows that this estimate is exact. Really, we have ($c$ is of [G]; $r$ of [G] is $q$ of the present paper\footnotemark \footnotetext{The equation of [G], 5.9.9 contains a misprint: the (2,1)-term of the first matrix is $1-c^{r+1}$ and not $1-c^{r}$ as it is printed.}) $\ord c=1$. After a change of the basis with the matrix $\left(\matrix c^{\frac{1+q+q^2}{1-q^2}}& 0\\0& c^{\frac{q}{1-q^2}} \endmatrix \right)$ the equation of [G], 5.9.9 becomes (2.4) of the present paper, with $\ord a_{12}=\frac{1}{q-1}, \ \ord a_{21}=-\frac{q^2}{q-1}$. We see that the example [G], 5.9.9 belongs to the initial point of the ray $\bar R_3$.
\medskip
{\bf Second proof of the equality $h^1(M(A))=h^1(M(A'))$ for $A$ of the form (3.4).}
\medskip
We consider a subset $U$ (described below) of $\Cal A_+$ on the $(v, \ t)$-coordinate plane, and we give a direct proof (without using Theorem 1.5) of
\medskip
{\bf Proposition 6.4.} Let $a_{12}, \ a_{21}$ be such that their $(v, \ t)$ belong to $U$, and let $A$ be of the form (3.4). In this case we have $h^1(M(A))=h^1(M(A'))$.
\medskip
{\bf Remark.} The restriction $(v, \ t)\in U$ is not essential, the same methods will give us a proof for any values of $(v, \ t)$.
\medskip
{\bf Proof.} Conditions $\ord a_{21}=v$, $\ord a_2=t$ mean that $\exists \ x, \ y\in \p$ such that $\ord x = \ord y=0$ and such that
$$a_{21}=\th^{-v}x, \ \ \ a_2=\th^{-t}y$$
We have $$a_{12}=(\frac{\th^q}{a_{21}}+\frac{\th^{q^2}}{a_{21}^{q^2}}-a_2)^{\frac1q}=\th^{1+\frac1q v}x^{-\frac1q}+\th^{qv+q}x^{-q}-\th^{-\frac1q t}y^{\frac1q}$$ For $(v,t)\in \Cal A_+$ we have $-(1+\frac1q v)>-(qv+q)<\frac1q t$, hence $u=\ord a_{12}=-(qv+q)$. We have $(v, \ t)\in \Cal A_+ \implies (u,v) \in \bar R_3$.
\medskip
We denote objects for the dual t-motive by prime, i.e. $a_{12}'=a_{21}$, $a_{21}'=a_{12}$ etc. We have\footnotemark \footnotetext{The reader can see that in this elementary calculation some terms are eliminated, similarly as teeth of a gear enter one into another. This is typical for mathematics; this shows symmetry of the present construction which is not seen in (3.9), (3.10) of [GL18].}
$$a_2'=\frac{\th^q}{a_{21}'}+\frac{\th^{q^2}}{{a_{21}'}^{q^2}}-{a_{12}'}^q=\frac{\th^q-{a_{12}'}^qa_{21}'}{a_{21}'}+\frac{\th^{q^2}}{{a_{21}'}^{q^2}}= $$ $$ =\frac{\th^{(-q+\frac1q)v+1}x^{q-\frac1q}+\th^{-qv-\frac1q t}x^qy^{\frac1q}}{a_{21}'}+\frac{\th^{q^2}}{{a_{21}'}^{q^2}}$$
We consider only the case $t<-v-q$. We denote $\Cal A_+ \cap \{t<-v-q\}$ by $U$. For $(v,t)\in U$ we have
$$\ord \th^{-qv-\frac1q t}x^qy^{\frac1q} < \ord \th^{(-q+\frac1q)v+1}x^{q-\frac1q}$$ hence $$t'=\ord a_2'= \ord \frac{\th^{-qv-\frac1q t}x^qy^{\frac1q}}{a_{21}'}= 2qv+\frac1q t+q$$ (because $\ord \frac{\th^{q^2}}{{a_{21}'}^{q^2}} > \ord \frac{\th^{-qv-\frac1q t}x^qy^{\frac1q}}{a_{21}'}$).
\medskip
We get that for $(v,t)\in U$ the numbers $v', \ t'$ are given by the formulas
$$v'=-qv-q$$
$$t'=2qv+\frac1q t+q$$
The easiest way to find the image of $U$ under this linear transformation is to find images of some points $P_*$ on borders of $U$. We have:
\medskip
$P_1$: $v=q, \ t=-q^3$ (side of both $U$, $\Cal A_+$). Its image: $v'=-q^2-q, \ t'=q^2+q$ (side of $\Cal A_-$).
\medskip
$P_2$: $v=\frac1{q-1}, \ t=-\frac{q^2}{q-1}$ (vertex of both $U$, $\Cal A_+$). Its image: $v'=-\frac{q^2}{q-1}, \ t'=\frac{q^2}{q-1}$ (vertex of $\Cal A_-$).
\medskip
$P_3$: $v=\frac1{q-1}, \ t=\frac{-q^2+q-1}{q-1}$ (vertex of $U$, side of $\Cal A_+$). Its image: $v'=-\frac{q^2}{q-1}, \ t'=\frac{q^3+q-1}{q^2-q}$ (side of $\Cal A_-$).
\medskip
$P_4$: $v=q, \ t=-2q$ (side of $U$, interior of $\Cal A_+$). Its image: $v'=-q^2-q, \ t'=2q^2+q-2$ (interior of $\Cal A_-$).
\medskip
We see that the image of $U$ under this linear transformation is in $\Cal A_-$, as it must be. $\square$
\medskip
{\bf 7. Case of $A$ of the form (3.5).}
\medskip
For simplicity, we consider the case $q=2$. We use notations $u=\ord a_{12}$, $v=\ord a_{21}$ as in Section 6. Since $a_{22}=0$ the value of $a_2$ for the present case is the same as in Section 6, hence we have the same $R_i$, $D_{ij}$ ($i,j=1,2,3)$ as in (6.0). We have: $\ord a_4=-4v-12$; \ \ $\ord a_3=-4v-8$; \ \ $\ord a_2\ge \min(-4v-4; -v-2; 2u)$; \ \ $\ord a_1=-v-1$; \ \ $\ord a_0=-v$. The tail coefficients are:
$\ord b_{14}=-4v-8$, \ \ $\ord b_{13}=-4v-4$, \ \ $\ord b_{12}=-4v$ for $v>0$, $\ord b_{12}=-v$ for $v<0$, $\ord b_{12}\ge0$ for $v=0$, $\ord b_{24}=-4v$.
\medskip
Case of Domain $D_{13}$.
\medskip
For $(u,v)\in D_{13}$ we have $\ord a_2=-4v-4$. If $v>-\frac13$ then the Newton polygon of (3.2.0) has vertices
\medskip
$(1,-v)$, $(4, -4v-4)$, $(8, -4v-8)$, $(16, -4v-12)$
\medskip
and $\ord x_{j0}$ are respectively $v+\frac43, v+\frac43, 1, \frac12$.
\medskip
If $v\le-\frac13$ then the Newton polygon of (3.2.0) has vertices
\medskip
$(1,-v)$, $(2, -v-1)$, $(8, -4v-8)$, $(16, -4v-12)$
\medskip
and $\ord x_{j0}$ are respectively $1, \frac12v+\frac76, \frac12v+\frac76, \frac12$.
\medskip
{\bf Proposition 7.1.} For the subdomain $(u,v)\in$ \{Domain $D_{13}$, $v\ge0$\}, we have $h^1=4$.
\medskip
{\bf Proof.} Let $x_0$ be a solution to (3.2.0). We consider its minimal solution $x_0, \ x_1, \ x_2,\dots$. We shall prove by induction that $\ord x_i\ge \frac12 i+\frac12$; this implies the proposition. For $i=0$ this is true: $\ord x_0\ge\frac12$. Let us prove that if this is true for some $i$ then this is true for $i+1$. First, we consider the case $i=0$. Ord's of the tail terms are:
\medskip
$\ord b_{14}x_{0}^{16}\ge -4v$; \ \ \ $\ord b_{13}x_{0}^{8}\ge -4v$; \ \ \ $\ord b_{12}x_{0}^{4}\ge -4v+2$,
\medskip
hence $g_1\ge -4v$. The Newton polygon has a point $(4, -4v-4)$, hence we get: $\ord x_1\ge1$, the induction supposition is true for $i=1$. So, now we consider the case $i>0$. For this case ord's of the tail terms are:
\medskip
$\ord b_{14}x_{i}^{16}\ge -4v+8i$;
\medskip
$\ord b_{13}x_{i}^{8}\ge -4v+4i$;
\medskip
$\ord b_{12}x_{i}^{4}\ge -4v+2+2i$,
\medskip
$\ord b_{24}x_{i-1}^{16}\ge -4v+8i$
\medskip
hence $g_{i+1}\ge -4v+2+2i$. Again, joining the points $(0,-4v+2+2i)$ and $(4, -4v-4)$ we get $\ord x_{i+1}\ge\frac12 (i+1)+\frac12$ --- the induction supposition for $i+1$. $\square$
\medskip
{\bf Proposition 7.2.} For the subdomain $(u,v)\in$ \{Domain $D_{13}$, $-\frac23< v\le0$\}, we have $h^1=4$.
\medskip
{\bf Proof.} Similar to the proof of Proposition 7.1, we use the same notations. As above we have $\ord x_0\ge\frac12$. Let us evaluate $\ord x_1$.
We have
\medskip
$\ord b_{14}x_{0}^{16}\ge -4v$;
\medskip
$\ord b_{13}x_{0}^{8}\ge -4v$;
\medskip
$\ord b_{12}x_{0}^{4}\ge -v+2$.
\medskip
If $-\frac13\le v\le0$ then $-4v< -v+2$, hence $g_1\ge -4v$. As in Proposition 7.1, we get $\ord x_1\ge1$ (if $v<-\frac13$ then the point $(4, -4v-4)$ is above the Newton polygon. In this case $\ord x_1\ge1$ also holds).

Now the induction supposition is the following: $\ord x_i\ge 4^{i-1}$. This is true for $i=1$. For a fixed $i$ we have:
\medskip
$\ord b_{14}x_{i}^{16}\ge -4v-8+4^{i+1}$;
\medskip
$\ord b_{13}x_{i}^{8}\ge -4v-4+2\cdot4^i$;
\medskip
$\ord b_{12}x_{i}^{4}\ge -v+4^i$,
\medskip
$\ord b_{24}x_{i-1}^{16}\ge -4v+4^i$
\medskip
For $i\ge1$ all these numbers are $\ge -v+4^i$. Joining the points $(0,-v+4^i)$ and $(1,-v)$ we get $\ord x_{i+1}\ge 4^i$. $\square$
\medskip
Propositions 7.1, 7.2 show that for all $(u,v)\in D_{13}$ we have $h^1=4$.
\medskip
Now we consider the case of Domain $D_{12}$.
For $(u,v)\in D_{12}$ we have $\ord a_2=-v-2$. First, we consider the case $v\ge-\frac43$.
\medskip
{\bf Proposition 7.3.} If $(u,v)\in D_{12}$ and $v\ge-\frac43$ then $h^1=4$.
\medskip
{\bf Proof.} If $(u,v)\in D_{12}$ and $v\ge-\frac43$ then the Newton polygon of (3.2.0) has vertices
\medskip
$(1,-v)$, $(2,-v-1)$, $(8, -4v-8)$, $(16, -4v-12)$
\medskip
and $\ord x_{j0}$ are respectively $1, \frac12v+\frac76, \frac12v+\frac76, \frac12$.
\medskip
We consider consecutively minimal solutions generated by $x_{j0}$ for all $j=1,\dots,4$.
\medskip
{\bf Lemma 7.3.1.} All minimal solutions generated by $x_{10}$ are small.
\medskip
{\bf Proof.} Since $\ord x_{10}=1$ we have
$\ord b_{14}x_{10}^{16}=-4v+8$, \ $\ord b_{13}x_{10}^{8}=-4v+4$, \ $\ord b_{12}x_{10}^{4}=-v+4$ and
\medskip
$\ord (b_{14}x_{10}^{16}+b_{13}x_{10}^{8}+b_{12}x_{10}^{4})= \min\{-4v+8,\ -4v+4, \ -v+4\}=-v+4$. Hence, $g_1=-v+4$ and $\ord x_{11}=v+(-v+4)=4$. By induction, continuing this calculation we get a small solution $x_{1i}$, $i\in \n N$, with $\ord x_{1,i+1}=4\cdot \ord x_{1i}$. Really, for $i=0$ we just showed that $\ord x_{10}=1$ and $\ord x_{11}=4$.
\medskip
Induction step $i \implies i+1$: the free term in the Newton polygon for $\ord x_{1,i+1}$ has the form $\ord b_{14}x_{1i}^{16}+b_{13}x_{1i}^{8}+b_{12}x_{1i}^{4}+b_{24}x_{1,i-1}^{16}$. Therefore its order is
\medskip
$\min \{-4v-8+16\cdot \ord x_{1i}, \ -4v-4+8\cdot \ord x_{1i}, \ -v+4\cdot \ord x_{1i}, \ -4v+16\cdot \ord x_{1,i-1}\}$ (because there exists only one minimal value, see below).
\medskip
By induction, the assumption $v\le-\frac23$ and the fact $\ord x_{1i}\ge0$, this is equal
$\min \{-4v-8+16\cdot \ord x_{1i}, \ -4v-4+8\cdot \ord x_{1i}, \ -v+4\cdot \ord x_{1i}, \ -4v+4\cdot \ord x_{1i}\}=-v+4\cdot \ord x_{1i}$ (the only minimal value).
\medskip
The leftmost segment of the Newton polygon gives us a solution $x_{1,i+1}$ with $\ord x_{1,i+1}=v+(-v+4\cdot \ord x_{1i})=4\cdot \ord x_{1i}$. $\square$
\medskip
{\bf Lemma 7.3.2.} All minimal solutions generated by $x_{j0}$, $j=2, \ 3$, are small.
\medskip
{\bf Proof.} Since $\ord x_{j0}=\frac12v+\frac76$ (here and below $j=2,\ 3$) we have
$\ord b_{14}x_{j0}^{16}=4v+\frac{32}3$, \ $\ord b_{13}x_{j0}^{8}=\frac{16}3$, \ $\ord b_{12}x_{j0}^{4}=v+\frac{14}3$ hence $g_1=v+\frac{14}3$. The leftmost segment of the Newton polygon gives us a solution $x_{j1}$ with $\ord x_{j1}=v+(v+\frac{14}3)=2v+\frac{14}3$ which is again $4\cdot \ord x_{j0}$.

As earlier, by induction, continuing this calculation we get small solutions $x_{ji}$, $i\in \n N$, with $\ord x_{j,i+1}=4\cdot \ord x_{ji}$. Really, for $i=0$ we just showed that $\ord x_{j1}=4\cdot \ord x_{j0}$.
\medskip
Induction step $i \implies i+1$: the free term in the Newton polygon for $\ord x_{j,i+1}$ has the form $\ord b_{14}x_{ji}^{16}+b_{13}x_{ji}^{8}+b_{12}x_{ji}^{4}+b_{24}x_{j,i-1}^{16}$. Therefore its order is
\medskip
$\min \{-4v-8+16\cdot \ord x_{ji}, \ -4v-4+8\cdot \ord x_{ji}, \ -v+4\cdot \ord x_{ji}, \ -4v+16\cdot \ord x_{j,i-1}\}$ (because there exists only one minimal value, see below).
\medskip
By induction, the assumption $-\frac43\le v\le-\frac23$ and the fact $\ord x_{ji}\ge0$, this is equal
$\min \{-4v-8+16\cdot \ord x_{ji}, \ -4v-4+8\cdot \ord x_{ji}, \ -v+4\cdot \ord x_{ji}, \ -4v+4\cdot \ord x_{ji}\}=-v+4\cdot \ord x_{ji}$ (the only minimal value).
\medskip
The leftmost segment of the Newton polygon gives us a solution $x_{j,i+1}$ with $\ord x_{j,i+1}=v+(-v+4\cdot \ord x_{ji})=4\cdot \ord x_{ji}$. $\square$
\medskip
{\bf Lemma 7.3.3.} All minimal solutions generated by $x_{40}$ are small.
\medskip
{\bf Proof.} It is completely analogous to the proofs of the above lemmas. Since $\ord x_{40}=\frac12$ we have
$\ord b_{14}x_{40}^{16}=-4v$, \ $\ord b_{13}x_{40}^{8}=-4v$, \ $\ord b_{12}x_{40}^{4}=-v+2$.

Hence $g_1=-v+2$, because $-v+2\le -4v$ if and only if $v\le -\frac23$. The leftmost segment of the Newton polygon gives us a solution $x_{41}$ with $\ord x_{41}=v+(-v+2)=2$ which is as before $4\cdot \ord x_{40}$.

As earlier, by induction, continuing this calculation we get small solutions $x_{4i}$, $i\in \n N$, with $\ord x_{4,i+1}=4\cdot \ord x_{4i}$. Really, for $i=0$ we just showed that $\ord x_{41}=4\cdot \ord x_{40}$.
\medskip
Induction step $i \implies i+1$: the free term in the Newton polygon for $\ord x_{4,i+1}$ has the form $\ord b_{14}x_{4i}^{16}+b_{13}x_{4i}^{8}+b_{12}x_{4i}^{4}+b_{24}x_{4,i-1}^{16}$. Therefore its order is
\medskip
$\min \{-4v-8+16\cdot \ord x_{4i}, \ -4v-4+8\cdot \ord x_{4i}, \ -v+4\cdot \ord x_{4i}, \ -4v+16\cdot \ord x_{4,i-1}\}$ (because there exists only one minimal value, see below).
\medskip
By induction, the assumption $-\frac43\le v\le-\frac23$ and the fact $\ord x_{4i}\ge\frac12$, this is equal
$\min \{-4v-8+16\cdot \ord x_{4i}, \ -4v-4+8\cdot \ord x_{4i}, \ -v+4\cdot \ord x_{4i}, \ -4v+4\cdot \ord x_{4i}\}=-v+4\cdot \ord x_{4i}$ (the only minimal value).
\medskip
The leftmost segment of the Newton polygon gives us a solution $x_{4,i+1}$ with $\ord x_{4,i+1}=v+(-v+4\cdot \ord x_{4i})=4\cdot \ord x_{4i}$. $\square$
\medskip
These 3 lemmas show that there are 4 linearly independent small solutions, i.e. $h^1=4$. $\square$
\medskip
{\bf Proposition 7.4.} If $(u,v)\in D_{12}$ and $v\le-\frac43$ then $h^1=4$.
\medskip
{\bf Proof.} If $(u,v)\in D_{12}$ and $v\le-\frac43$ then the Newton polygon of (3.2.0) has vertices
\medskip
$(1,-v)$, $(2,-v-1)$, $(4, -v-2)$, $(16, -4v-12)$
\medskip
and $\ord x_{j0}$ are respectively $1, \frac12, \frac14v+\frac56, \frac14v+\frac56$.
\medskip
For $j=1$, resp. 2, we have the same values of $\ord x_{ji}$ (minimal solutions) as in Lemmas 7.3.1, resp. 7.3.2, i.e. these minimal solutions are small.
\medskip
Let us consider the case $j=3,\ 4$, i.e. $\ord x_{j0}=\frac14v+\frac56$. We have
$\ord b_{14}x_{j0}^{16}=\frac{16}3$, \ $\ord b_{13}x_{j0}^{8}=-2v+\frac83$, \ $\ord b_{12}x_{j0}^{4}=\frac{10}3$. Since $v\le-\frac23$ the minimum of these three values is $\frac{10}3$ and we have $g_1=\frac{10}3$.
\medskip
Therefore, there are three cases of the Newton polygon:
\medskip
(a) If $v+\frac{10}3>1$ equivalently $v>-\frac73$ or $\ord x_{j0}>\frac14$, then the Newton polygon for $x_{j1}$ has the first two vertices $(0, \frac{10}3)$ and $(1,-v)$ and the order of the small solution $x_{j1}$ is $v+\frac{10}3=4\cdot \ord x_{j0}$.
\medskip
(b) If $1\ge v+\frac{10}3$ and $\frac{v+1+\frac{10}3}2>\frac12$, equivalently $-\frac73 \ge v>-\frac{10}3$ or $\frac14\ge \ord x_{j0}>0$, then the Newton polygon for $x_{j1}$ has the first two vertices $(0, \frac{10}3)$ and $(2,-v-1)$ and the order of the small solution $x_{j1}$ is $\frac12(v+\frac{10}3)=\frac12(4\cdot \ord x_{j0})+\frac12$.
\medskip
(c) If $\frac{v+1+\frac{10}3}2\le\frac12$, equivalently $v\le-\frac{10}3$ or $\ord x_{j0}\le0$, then the Newton polygon for $x_{j1}$ has the first two vertices $(0, \frac{10}3)$ and $(4,-v-2)$ and the order of the small solution $x_{j1}$ is $\frac14(v+\frac{10}3)+\frac12=\frac14(4\cdot \ord x_{j0})+\frac12=\ord x_{j0}+\frac12$.
\medskip
{\bf Lemma 7.4.1.} If there exists $i\in \n N$ such that $\ord x_{ji}=4\cdot \ord x_{j,i-1}>1$ then $\ord x_{j,i+1}=4\cdot \ord x_{ji}$.
\medskip
{\bf Proof.} By hypothesis we have $\ord x_{ji}=4\cdot \ord x_{j,i-1}$. Under this condition the ord's of the tail terms:

$-4v-8 + 16\cdot \ord x_{ji}, \ -4v-4+8\cdot \ord x_{ji}, \ -v+4\cdot \ord x_{ji}, \ -4v+4\cdot \ord x_{ji}$ are different, and conditions $v\le-\frac43$ and $\ord x_{ji}>1$ imply that their minimum is $-v+4\cdot \ord x_{ji}$. Since $-v-(-v+4\cdot \ord x_{ji})<-1$ then the Newton polygon for $x_{j,i+1}$ has the first two vertices $(0, -v+4\cdot \ord x_{ji})$ and $(1,-v)$ and the order of the minimal solution $x_{j,j+1}$ is $4\cdot \ord x_{ji}$. $\square$
\medskip
{\bf Lemma 7.4.2.} If there exists $i\in \n N$ such that $\ord x_{ji}=2\cdot \ord x_{j,i-1}+\frac12>\frac12$, then $\ord x_{j,i+1}=4\cdot \ord x_{ji}>2$.
\medskip
{\bf Proof.} The ord's of the tail terms are the following:

$-4v-8 + 16\cdot \ord x_{ji}, \ -4v-4+8\cdot \ord x_{ji}, \ -v+4\cdot \ord x_{ji}, \ -4v+16\cdot \ord x_{j,i-1}$.
By hypothesis we have $-4v-4+8\cdot \ord x_{ji}=-4v+16\cdot \ord x_{j,i-1}$ and with $\ord x_{ji}>\frac12$ we get
$-4v-4+8\cdot \ord x_{ji} < -4v-8 + 16\cdot \ord x_{ji}$. Because of $v\le-\frac43$ and $\ord x_{ji}>0$ we have $-v+4\cdot \ord x_{ji} < -4v-4+8\cdot \ord x_{ji}$. We get

$\min \{ -4v-8 + 16\cdot \ord x_{ji}, \ -4v-4+8\cdot \ord x_{ji}, \ -v+4\cdot \ord x_{ji}, \ -4v+16\cdot \ord x_{j,i-1} \} = -v+4\cdot \ord x_{ji}$ (the only minimal value). Since $-v-(-v+4\cdot \ord x_{ji})<-2$ then the Newton polygon for $x_{j,i+1}$ has the first two vertices $(0, -v+4\cdot \ord x_{ji})$ and $(1,-v)$ and the order of the minimal solution $x_{j,j+1}$ is $4\cdot \ord x_{ji}$ which is $\ge2$. $\square$
\medskip
{\bf Lemma 7.4.3.} If $\ord x_{ji}=\ord x_{j,i-1}+\frac12\le\frac12$ and if $\ord x_{ji}=\frac14v+c$ where $c>\frac23$ then

$$\ord x_{j,i+1}\ =\ \matrix 4\cdot \ord x_{ji} & \hbox{ if } & \ord x_{ji}>\frac14 \\
2\cdot \ord x_{ji}+\frac12 & \hbox{ if } & 0<\ord x_{ji}\le\frac14 \\
 \ord x_{ji}+\frac12 & \hbox{ if } & \ord x_{ji}\le0 \endmatrix $$

{\bf Proof.} The ord's of the tail terms are the following:

$-4v-8 + 16\cdot \ord x_{ji}, \ -4v-4+8\cdot \ord x_{ji}, \ -v+4\cdot \ord x_{ji}, \ -4v+16\cdot \ord x_{j,i-1}$.
By hypothesis we have $-4v-8+16\cdot \ord x_{ji}=-4v+16\cdot \ord x_{j,i-1}$ and with $\ord x_{ji}\le\frac12$ we get
$-4v-4+8\cdot \ord x_{ji} \ge -4v-8 + 16\cdot \ord x_{ji}$. Because of $\ord x_{ji}=\frac14v+c$ we have $-4v-8 + 16\cdot \ord x_{ji}=-8 + 16c$ and $-v+4\cdot \ord x_{ji}=4c$.
\medskip
Condition $c>\frac23$ implies $-4v-8 + 16\cdot \ord x_{ji}>-v+4\cdot \ord x_{ji}$. Therefore $\min \{ -4v-8 + 16\cdot \ord x_{ji}, \ -4v-4+8\cdot \ord x_{ji}, \ -v+4\cdot \ord x_{ji}, \ -4v+16\cdot \ord x_{j,i-1} \} = -v+4\cdot \ord x_{ji}$ (the only minimal value). If $-v-(-v+4\cdot \ord x_{ji})<-1$, equivalently $\ord x_{ji}>\frac14$, then the Newton polygon for $x_{j,i+1}$ has the first two vertices $(0, -v+4\cdot \ord x_{ji})$ and $(1,-v)$ and the order of the minimal solution $x_{j,j+1}$ is $4\cdot \ord x_{ji}$ which is $\ge2$.
\medskip
If $-v-(-v+4\cdot \ord x_{ji})\ge-1$ and $\frac{-v-1-(-v+4\cdot \ord x_{ji})}2<-\frac12$, equivalently $0<\ord x_{ji}\le\frac14$, then the Newton polygon for $x_{j,i+1}$ has the first two vertices $(0, -v+4\cdot \ord x_{ji})$ and $(2,-v-1)$ and the order of the minimal solution $x_{j,j+1}$ is $\frac{4\cdot \ord x_{ji}+1}2=2\cdot \ord x_{ji}+\frac12$.
\medskip
If $(-v-1-(-v+4\cdot \ord x_{ji})\ )/2\ge-\frac12$, equivalently $0\ge\ord x_{ji}$, then the Newton polygon for $x_{j,i+1}$ has the first two vertices $(0, -v+4\cdot \ord x_{ji})$ and $(4,-v-2)$ and the order of the minimal solution $x_{j,j+1}$ is $(4\cdot \ord x_{ji}+2)/4=\ord x_{ji}+\frac12$. $\square$
\medskip
Now we can finish the proof. We have $\ord x_{j0}=\frac14v+\frac56$, $j=3, \ 4$. Condition $v\le-\frac43$ implies $\ord x_{j0}\le\frac12$. Lemma 7.4.3 implies: If $\ord x_{j0}>\frac14$ then $\ord x_{j1}=4\cdot \ord x_{j0}>1$. Lemma 7.4.1 shows that the corresponding minimal solution is small, it satisfies $\ord x_{ji}=4^i\cdot \ord x_{j0}$.
\medskip
If $0<\ord x_{j0}\le\frac14$ then $\ord x_{j1}=2\cdot \ord x_{j0}+\frac12>\frac12$. Applying Lemma 7.4.2 we get $\ord x_{j2}=4\cdot \ord x_{j1}>2$. Lemma 7.4.1 shows that the corresponding minimal solution is small, it satisfies $\ord x_{ji}=4^{i-1}\cdot \ord x_{j1}$.
\medskip
If $\ord x_{j0}\le0$ then $\ord x_{j1}=\ord x_{j0}+\frac12$. There exists a number $k\in \n N$ such that $0<\ord x_{j0}+k\cdot\frac12\le\frac12$. By Lemma 7.4.3 we have $\ord x_{jk}=\ord x_{j,k-1}+\frac12$. Also by Lemma 7.4.3 if $\ord x_{jk}\le\frac14$ then $\ord x_{j,k+1}=2\cdot \ord x_{jk}+\frac12>\frac12$ and therefore $\ord x_{j,k+2}=4\cdot \ord x_{j,k+1}>1$. Hence, we get a small solution with $\ord x_{j,k+i}=4^{i-1}\cdot \ord x_{j,k+1}$. Otherwise, if $\ord x_{jk}>\frac14$ then $\ord x_{j,k+1}=4\cdot \ord x_{j,k}>1$ and therefore we get a small solution with $\ord x_{j,k+i}=4^{i}\cdot \ord x_{j,k}$. $\square$
\medskip
These propositions show that for all $(u,v)\in D_{12}$ we have $h^1=4$.
\medskip
{\bf Remark 7.4.4.} End of the proof of Proposition 7.4 (case of Lemma 7.4.3 and the very end of the proof) shows that two types of the Newton polygon and respectively two types of the growth of $\ord x_i$ can occur: linear for small $i$, $i<k$ and exponential for $i\ge k$.
\medskip
Now we consider the case $(u,v)\in D_{23}$. We get a result only for a subset of $D_{23}$ denoted in Section 9 by $D_{23}$1.1.2. It is defined by the inequalities $v\le -\frac43$, $-2u-4\ge v\ge -2u-3$.
\medskip
{\bf Proposition 7.5.} For all $(u,v)\in$ $D_{23}$1.1.2 we have $h^1=4$.
\medskip
{\bf Proof.} We have (see Section 9): If $v\le -\frac43$ then the convex hull of the points
\medskip
$(1,-v)$, $(2, -v-1)$, $(8, -4v-8)$, $(16, -4v-12)$
\medskip
consists of the points $(1,-v)$, $(2, -v-1)$, $(16, -4v-12)$ (case $D_{23}$1).
If $v\le -\frac43$, $2u \le -v-3$ the Newton polygon of (3.2.0) has vertices $(1,-v)$, $(4, 2u)$, $(16, -4v-12)$ (case $D_{23}$1.1 of Section 9).
\medskip
We have $\ord x_{j0}$ are respectively $\frac16u+\frac13v+1, \frac16u+\frac13v+1, -\frac23u-\frac13v, -\frac23u-\frac13v$.
\medskip
We prove by induction by $i$ that for all $(u,v)\in$ $D_{23}$1.1.2, for all solutions $\{x\}$, for all $i\ge0$ we have: $\ord x_{i}\ge (\frac16-\frac12i)u+(\frac13-\frac14i)v+1$. This is true for $i=0$. Using the induction assumption for $i$, we get:
\medskip
$\ord b_{14}x_{i}^{16}\ge(\frac83-8i)u+(\frac43-4i)v+8$,
\medskip
$\ord b_{13}x_{i}^{8}\ge(\frac43-4i)u+(-\frac43-2i)v+4$,
\medskip
$\ord b_{12}x_{i}^{4}\ge(\frac23-2i)u+(\frac13-i)v+4$,
\medskip
$\ord b_{24}x_{i-1}^{16}\ge(\frac{32}3-8i)u+(\frac{16}3-4i)v+16$.
\medskip
For all $(u,v)\in$ $D_{23}$1.1.2, for all $i\ge0$ a minimal of these four numbers is the third one (for $i=0$ the fourth number is not considered, and the first number is equal to the third one). This means that $g_{i+1}$ is $\ge (\frac23-2i)u+(\frac13-i)v+4$. Further, for all $(u,v)\in$ $D_{23}$1.1.2, for all $i\ge0$ we have $(\frac23-2i)u+(\frac13-i)v+4>\frac83u+\frac43v+4$, hence $\ord x_{i+1}\ge$ is the minus slope of the segment $(0, (\frac23-2i)u+(\frac13-i)v+4); (4,2u)$ which is equal $(\frac16-\frac12(i+1))u+(\frac13-\frac14(i+1))v+1$. This implies the induction supposition.
\medskip
Finally, we have $\forall \ (u,v)\in$ $D_{23}$1.1.2 $(\frac16-\frac12i)u+(\frac13-\frac14i)v+1$ tends to $+\infty$ as $i$ tends to $+\infty$, hence any $x_{0}$ generates a small solution, hence all Anderson t-motives in $D_{23}$1.1.2 are uniformizable. $\square$
\medskip
{\bf 8. Case of jump of valuation of the terms of $a_2$.}
\medskip
Here we continue to study the case of $A$ of the form (3.5), but now we consider the case when there exists a jump of valuation of the terms of $a_2$, i.e. ord's of some of the terms of $a_2$ are equal. This can occur if $(u,v)\in R_1\cup R_2\cup R_3$, where $R_i$ are from Section 6. As in Section 3, we denote $t:=\ord a_2$.
\medskip
{\bf Proposition 8.1.} If $(v,t)$ satisfy $0\ge v<2$, $t\ge -\frac{16}7 v -\frac{24}7$ (this is a subset of $D_{23}$3.3 of Section 9) then $h^1=4$.
\medskip
{\bf Proof.} Condition $v\ge-\frac13$, $t\ge -\frac{16}7 v -\frac{24}7$ means that the vertices of the convex hull of the Newton polygon of (3.2.0) are $(1,-v)$, $(8, -4v-8)$, $(16, -4v-12)$ and
Ord's of $x_{j0}$ are $\frac12$, $\frac37v+\frac87$ (three times). In particular, $\forall \ i=1,\dots,4$ we have $\ord x_{i0}\ge\frac12$. Before starting induction, we evaluate $\ord x_{i1}$, $\ord x_{i2}$. We have
\medskip
$\ord b_{14}x_{i0}^{16}\ge-4v$,
\medskip
$\ord b_{13}x_{i0}^{8}\ge-4v$,
\medskip
$\ord b_{12}x_{i0}^{4}\ge-4v+2$ (because $v\ge0$), hence $\ord g_1\ge-4v$. Joining the vertices $(0,\ord g_1)$ and $(8,-4v-8)$ we get $\ord x_{i1}\ge1$.
\medskip
Further, we have
\medskip
$\ord b_{14}x_{i1}^{16}\ge-4v+8$,
\medskip
$\ord b_{13}x_{i1}^{8}\ge-4v+4$,
\medskip
$\ord b_{12}x_{i1}^{4}\ge-4v+4$,
\medskip
$\ord b_{24}x_{i0}^{16}\ge-4v+8$. By the same reason, we get $\ord x_{i2}\ge\frac32$.

The first induction supposition is: $\ord x_{ij}\ge2-\frac1{2^{j-1}}$. This is true for $j=1,\ 2$. We assume that this is true for $j-1, \ j$, and prove that this is true for $j+1$. We have $\ord x_{i,j+1}\ge \min(\ (\ord b_{14}x_{ij}^{16}-(-4v-8))/8, \ (\ord b_{13}x_{ij}^{8}-(-4v-8))/8,$

$ \ (\ord b_{12}x_{ij}^{4}-(-4v-8))/8, \ (\ord b_{24}x_{i,j-1}^{16}-(-4v-8))/8\ )$. These four numbers satisfy inequalities
\medskip
$(\ord b_{14}x_{ij}^{16}-(-4v-8))/8 \ge2(2-\frac1{2^{j-1}})$,
\medskip
$(\ord b_{13}x_{ij}^{8}-(-4v-8))/8 \ge2-\frac1{2^{j-1}}+\frac12$,
\medskip
$(\ord b_{12}x_{ij}^{4}-(-4v-8))/8\ge2-\frac1{2^{j}}$,
\medskip
$(\ord b_{24}x_{i,j-1}^{16}-(-4v-8))/8\ge5-\frac1{2^{j-3}}$. The minimal is the third number; this implies the induction supposition.
\medskip
Now we consider the segment of the Newton polygon $(0,g_i)$ --- $(1,-v)$. We denote $\ve=2-v$. There exists $j_0$ such that $\ord x_{i,j_0-1}>2-\frac\ve4$. Let us evaluate $\ord x_{i,j_0+1}$, $\ord x_{i,j_0+2}$: for them the point $(1,-v)$ belongs to the convex hull of the Newton polygon for these $x_{ij}$. Really,

$\min (\ \ord b_{14}x_{i,j_0}^{16}, \ \ord b_{13}x_{i,j_0}^{8}, \ \ord b_{12}x_{i,j_0}^{4}, \ \ord b_{24}x_{i,,j_0-1}^{16})=\ord b_{12}x_{i,j_0}^{4}$

(because $\ord x_{i,j_0-1}\ge\frac32$), and the slope of the segment of the Newton polygon $(0,g_i)$ --- $(1,-v)$ is $\ge 2+2\ve$, i.e. $\ord x_{i,j_0+1}\ge 2+2\ve$. Analogously, $\ord x_{i,j_0+2}\ge 2+2\ve$ (we use the above arguments for $j_0+1$ instead of $j_0$).
\medskip
The second induction supposition is: $\ord x_{i,j_0+j}\ge 2+2^{j-1}\ve$ ($j>0$). Really, this is true for $j=1, \ 2$. We assume that this is true for $j-1, \ j$, and prove that this is true for $j+1$. We have
\medskip
$\ord b_{14}x_{i,j_0+j}^{16}-(-v) \ge-3v+24+2^{j-1}\cdot16\ve>2+2^j$,
\medskip
$\ord b_{13}x_{i,j_0+j}^{8}-(-v) \ge-3v+12+2^{j-1}\cdot8\ve>2+2^j$,
\medskip
$\ord b_{12}x_{i,j_0+j}^{4}-(-v)\ge-3v+8+2^{j-1}\cdot4\ve>2+2^j$,
\medskip
$\ord b_{24}x_{i,j_0+j-1}^{16}-(-v)\ge-3v+32+2^{j-2}\cdot16\ve>2+2^j$.

This proves the induction supposition, and hence the proposition. $\square$
\medskip
{\bf Proposition 8.2.} If $(v,t)$ satisfy $v\ge2$, $t\ge -\frac{16}7 v -\frac{24}7$ (this is a subset of $D_{23}$3.3) then $h^1=0$.
\medskip
{\bf Proof.} First, we consider the maximal value of $\ord x_0$ and their minimal solution $x_i$. We denote $\al_i:=\ord x_i-2$. We have $\al_0=\frac37v-\frac67\ge0$. The induction supposition is $\al_{i+1}=\al_i/2$. Let us prove it. We have

\medskip
$\ord b_{14}x_{i}^{16}=-4v+16\al_i+24$,
\medskip
$\ord b_{13}x_{i}^{8}=-4v+8\al_i+12$,
\medskip
$\ord b_{12}x_{i}^{4}=-4v+4\al_i+8$,
\medskip
$\ord b_{24}x_{i-1}^{16}=-4v+32\al_i+16$.
\medskip

According the induction supposition, we have $\al_i>0$, hence the above numbers are different and $g_{i+1}=-4v+4\al_i+8$. The convex hull of $(0,g_{i+1}), \ (1,-v), \ (2,-v-1), \ (4,t), \ (8,-4v-8), \ (16,-4v-12)$ consists of the points $(0,g_{i+1}), \ (8,-4v-8), \ (16,-4v-12)$ --- this follows immediately from $\frac37v-\frac67\ge\al_i\ge0$ and $v\ge2$. Hence, the maximal value of $\ord x_{i+1}$ is $2+\al_i/2$, hence the supposition.
\medskip
This means that these minimal solutions are not small.
\medskip
Let us consider the case of $\ord x_0=\frac12$ and its minimal solution. We have in this case
\medskip
$\ord b_{14}x_{0}^{16}=-4v$,
\medskip
$\ord b_{13}x_{0}^{8}=-4v$,
\medskip
$\ord b_{12}x_{0}^{4}=-4v+2$, i.e. a jump can occur. Hence, we have to consider one more term in an approximation. We denote $c_2:=a_2\th^{-\frac{16}7 v -\frac{24}7}$, i.e. $\ord c_2\ge0$. The first approximation to $x_0$ is $x_0\approx \frac{a_{11}^{1/2}}{\th}$ (it is obtained if we consider only the first two terms of (3.2.0)). Hence, we let $\De_0:=x_0-\frac{a_{11}^{1/2}}{\th}$, i.e. $x_0=\De_0+\frac{a_{11}^{1/2}}{\th}$. Substituting this value to (3.2.0) we get$$a_4\De_0^{16}+a_{3}\De_0^{8}+a_2\De_0^{4}+a_1\De_0^2+a_0\De_0 + (\th^{\frac{16}7 v -\frac{4}7}a_{11}^{2}c_2+\frac{a_{11}^{2}}{a_{21}\th^2} +\frac{a_{11}^{1/2}}{a_{21}\th})=0\eqno{(8.2.1)}$$
We have: $\ord \th^{\frac{16}7 v -\frac{4}7}a_{11}^{2}c_2\ge \frac{-16v-10}7$ which is always less than $\ord \frac{a_{11}^{2}}{a_{21}\th^2}$, $\ord \frac{a_{11}^{1/2}}{a_{21}\th}$, hence ord of the free term of (8.2.1) is $\ge \frac{-16v-10}7$. Comparing with a Newton point $(8,-4v-8)$ we get that $\ord \De_0\ge \frac{6 v +23}{28}\ge\frac54$, because $v\ge2$.
\medskip
Now we substitute the obtained value of $x_0$ to (3.2.1). The sum of the tail terms is
$$\frac{\th^8+\th^4}{a_{21}^4}(\frac{a_{11}^8}{\th^{16}}+\De_0^{16})+\frac{a_{11}^4}{a_{21}^4}(\frac{a_{11}^4}{\th^{8}}+\De_0^{8})+(\frac1{a_{21}} +\frac1{a_{21}^4})(\frac{a_{11}^2}{\th^{4}}+\De_0^{4})\eqno{(8.2.2)}$$ Taking into consideration that $\ord \De_0\ge \frac54$ we get that the term $\frac1{a_{21}^4}\frac{a_{11}^2}{\th^{4}}$ of (8.2.2) has the minimal ord $=-4v+2$, all other terms of (8.2.2) have the higher ord's. This means that $g_1=-4v+2$ and $\ord x_1=\frac54$.
\medskip
The continuation is simpler because there will be no more jumps of valuation. The induction supposition is: for $i\ge1$ we have $\ord x_i=2-\frac3{2^{i+1}}$. This is true for $i=1$. Let us assume that this is true for some $i$, and prove that this is true for $i+1$. We have
\medskip
$\ord b_{14}x_{i}^{16}=-4v+24-\frac3{2^{i-3}}$,
\medskip
$\ord b_{13}x_{i}^{8}=-4v+12-\frac3{2^{i-2}}$,
\medskip
$\ord b_{12}x_{i}^{4}=-4v+8-\frac3{2^{i-1}}$
\medskip
$\ord b_{24}x_{0}^{16}=-4v+8$ for $i=1$, $\ord b_{24}x_{i-1}^{16}=-4v+24+\frac3{2^{i-2}}$ for $i>1$.
\medskip
For $i\ge1$ the minimal of the above four numbers is $-4v+8-\frac3{2^{i-1}}$. The leftmost segment of the convex hull of $(0, -4v+8-\frac3{2^{i-1}})$, $(1,-v)$, $(2, -v-1)$, $(8, -4v-8)$, $(16, -4v-12)$ is [$(0, -4v+8-\frac3{2^{i-1}})$, $(8, -4v-8)$], hence $\ord x_{i+1}=2-\frac3{2^{i}}$ --- the induction supposition is proved. We see that a minimal solution is not small.
\medskip
Let us prove that there is no small solutions at all (in principle it can happen that a linear combination of non-small solutions is small). The proof is exactly the same as the proof of [GL18], Lemma 4.6, so we repeat it here. Namely, we denote by $\{x_j\}$ $(j=1,2,3)$ linearly independent minimal solutions corresponding to $\ord x_{j0}=\frac37v+\frac87$, and by $\{x_4\}$ a minimal solution corresponding to $\ord x_{40}=\frac12$. Let us assume that $\exists \ C_1,...,C_4\in \n F_2[[T]]$ such that $\sum_{j=1}^4 C_j \{x_j\}$ is a small solution. We consider $S_{123}:=\sum_{j=1}^3 C_j \{x_j\}$, we denote $S_{123}=\sum_{i=0}^\infty \bar x_{1,2,3;i}T^i$. We have: $\ord \bar x_{1,2,3;i}\ge2$, because $\forall \ i$ elements $\bar x_{1,2,3;i}$ are linear combinations of $x_{jk}$ for $j=1,2,3$, $k\le i$ with coefficients in $\n F_2$.

Further, we denote $S_{4}:=C_4 \{x_4\}=\sum_{i=0}^\infty \bar x_{4i}T^i$. The above considerations show that $\forall \ i$ $\ord x_{4i}$ are different and $\frac12\le \ord x_{4i}<2$, hence $\forall \ i$ we have $\frac12\le \ord \bar x_{4i}<2$. This means that $\sum_{j=1}^4 C_j \{x_j\}=S_{123}+S_4$ cannot be a small solution. $\square$
\medskip
{\bf Proposition 8.3.} If $(v,t)$ satisfy conditions $-\frac43<v<-\frac13$, $t\ge -4v-6$  then $h^1=4$.
\medskip
{\bf Remark.} The conditions $-\frac43<v<-\frac13$, $t\ge -4v-6$ correspond to the types $D_{23}$2a.3, $D_{23}$2b.2, $D_{23}$2b.3, $D_{23}$2.5 of Section 9.
\medskip
{\bf Proof.} For all these cases we have $\ord x_{i0}\ge\frac12$. First, we evaluate $\ord x_{i1}$. We have
\medskip
$\ord b_{14}x_{i0}^{16}\ge-4v$,
\medskip
$\ord b_{13}x_{i0}^{8}\ge-4v$,
\medskip
$\ord b_{12}x_{i0}^{4}\ge-v+2$, hence $\ord x_{i1}\ge\min(-3v,2)=1$. Induction supposition for any $j$: $ \ord x_{ij}\ge2^{j-1}$. This is true for $j=0, \ 1$; let us prove that it holds for $j+1$. We have
\medskip
$\ord b_{14}x_{ij}^{16}\ge\frac43-8+16\cdot2^{j-1}$,
\medskip
$\ord b_{13}x_{ij}^{8}\ge\frac43-4+8\cdot2^{j-1}$,
\medskip
$\ord b_{12}x_{ij}^{4}\ge\frac13+4\cdot2^{j-1}$
\medskip
$\ord b_{24}x_{i,j-1}^{16}\ge\frac43+16\cdot2^{j-2}$
\medskip
To prove the induction supposition, we must prove that for $j\ge1$ we have: $\min(\frac43-8+16\cdot2^{j-1},\ \frac43-4+8\cdot2^{j-1}, \ \frac13+4\cdot2^{j-1}, \frac43+16\cdot2^{j-2})-\frac43\ge 2^j$. This is straightforward. $\square$
\medskip
{\bf 9. Complementary information for the case of $A$ from (3.5).}
\medskip
Here we give some information and results of preliminary calculations that can be useful for further research.
\medskip
Case Domain $D_{23}$.
\medskip
Let us consider first the convex hull of the points
\medskip
$(1,-v)$, $(2, -v-1)$, $(8, -4v-8)$, $(16, -4v-12)$.
\medskip
{\bf Case }$D_{23}$1. If $v\le -\frac43$ it consists of the points $(1,-v)$, $(2, -v-1)$, $(16, -4v-12)$. Depending on $u$ we have 3 possibilities:
\medskip
If $v\le -\frac43$, $2u \le -v-3$ the Newton polygon has vertices $(1,-v)$, $(4, 2u)$, $(16, -4v-12)$
\medskip
and $\ord x_{j0}$ are respectively $\frac16u+\frac13v+1, \frac16u+\frac13v+1, -\frac23u-\frac13v, -\frac23u-\frac13v$.
\medskip
If $v\le -\frac43$, $-\frac{10}7v-\frac{18}7\ge 2u \ge -v-3$ the Newton polygon has vertices $(1,-v)$, $(2, -v-1)$, $(4, 2u)$, $(16, -4v-12)$
\medskip
and $\ord x_{j0}$ are respectively $\frac16u+\frac13v+1, \frac16u+\frac13v+1, -u-\frac12v-\frac12, 1$.
\medskip
If $v\le -\frac43$, $-\frac{10}7v-\frac{18}7\le 2u$ the Newton polygon has vertices $(1,-v)$, $(2, -v-1)$, $(16, -4v-12)$
\medskip
and $\ord x_{j0}$ are respectively $\frac3{14}v+\frac{11}{14}, \frac3{14}v+\frac{11}{14}, \frac3{14}v+\frac{11}{14},  1$.
\medskip
In Domain $D_{23}$, $v<-\frac43$ we have always $-\frac{10}7v-\frac{18}7\ge 2u$, hence really we have 2 possibilities:
\medskip
Case $D_{23}$1.1. Domain: $v\le -\frac43$, $v \le -2u-3$.
\medskip
Case $D_{23}$1.2. Domain: $v\le -\frac43$, $-2u-2\ge v \ge -2u-3$.
\medskip
If we consider $t=\ord a_2$ we get the third possibility:
\medskip
Case $D_{23}$1.3. $v\le-\frac43$, $t\ge -\frac{10}7 v -\frac{18}7$.
\medskip
The vertices of the Newton polygon for these cases are:
\settabs 6 \columns
\medskip
\+Case $D_{23}$1.1&$(1,-v)$& & $(4, 2u)$& &$(16, -4v-12)$\cr
\medskip
\+Case $D_{23}$1.2&$(1,-v)$& $(2, -v-1)$ & $(4, 2u)$& & $(16, -4v-12)$\cr
\medskip
\+Case $D_{23}$1.3&$(1,-v)$& $(2,-v-1)$& && $(16, -4v-12)$\cr
\medskip
Ord's of $x_{j0}$:
\settabs 5 \columns
\medskip
\+Case $D_{23}$1.1&$\frac16u+\frac13v+1$&$\frac16u+\frac13v+1$&$ -\frac23u-\frac13v$&$ -\frac23u-\frac13v$\cr
\medskip
\+Case $D_{23}$1.2& $\frac16u+\frac13v+1$&$  \frac16u+\frac13v+1$&$  -u-\frac12v-\frac12$&$  1$\cr
\medskip
\+Case $D_{23}$1.3& $\frac3{14}v+\frac{11}{14}$&$\frac3{14}v+\frac{11}{14}$&$\frac3{14}v+\frac{11}{14}$&$  1$\cr
\medskip
Here we start a calculation of ord $x_1$ for some cases.
\medskip
Case $D_{23}$1.1, $j=1$. We have $\ord x_{10}=\frac16u+\frac13v+1$,
\medskip
$\ord b_{14}x_{10}^{16}=\frac83u+\frac43v+8$,
\medskip
$\ord b_{13}x_{10}^{8}=\frac43u-\frac43v+4$,
\medskip
$\ord b_{12}x_{10}^{4}=\frac23u+\frac13v+4$.
\medskip
In all points of $D_{23}$1.1 we have $\frac83u+\frac43v+8<\frac43u-\frac43v+4$,
hence, the domain $D_{23}$1.1 consists of two subdomains $D_{23}$1.1.1 ($v<-2u-4$) and $D_{23}$1.1.2 ($v>-2u-4$). We have the $g_1=\frac83u+\frac43v+8$ in $D_{23}$1.1.1, $g_1=\frac23u+\frac13v+4$ in $D_{23}$1.1.2.
\medskip
Case $D_{23}$1.1.1. The line defined by points $(4, 2u)$, $(16, -4v-12)$ (vertices of Newton polygon) crosses the $v$-axis at the point $(0, \frac83u+\frac43v+4)$. Further, the segment $[(0,\frac83u+\frac43v+8), \ (4, 2u)]$ crosses the line $u=1$ at the point $(1,\frac52u+v+6)$ which is under the point $(1,-v)$ for all $(u,v)\in $ $D_{23}$1.1.1. Hence, the segment $[(0,\frac83u+\frac43v+8), \ (4, 2u)]$ is a part of the Newton polygon of $x_1$, and $\ord x_{11}=\frac16u+\frac13v+2$ for all points of $D_{23}$1.1.1.
\medskip
Let us continue for $\ord x_{12}$:
\medskip
$\ord b_{14}x_{11}^{16}=\frac83u+\frac43v+24$,
\medskip
$\ord b_{13}x_{11}^{8}=\frac43u-\frac43v+12$,
\medskip
$\ord b_{12}x_{11}^{4}=\frac23u+\frac13v+8$,
\medskip
$\ord b_{24}x_{10}^{16}=\frac83u+\frac43v+16$.
\medskip
In all points of $D_{23}$1.1 we have $\frac83u+\frac43v+16<\frac43u-\frac43v+12$, $\frac83u+\frac43v+16<\frac83u+\frac43v+24$, so we need to compare $\frac83u+\frac43v+16$ and $\frac23u+\frac13v+8$. We did not finish this calculation.
\medskip
{\bf Case }$D_{23}$2. If $-\frac13 \ge v\ge -\frac43$ the convex hull of $(1,-v)$, $(2, -v-1)$, $(8, -4v-8)$, $(16, -4v-12)$ consists of all these points. Depending on $u$ we have 5 domains on $(u,v)$-plane defined as follows:
\medskip
Case $D_{23}$2.1. $-\frac13 \ge v\ge -\frac43$, $v\le -\frac12u-\frac32$, $v\le -2u-3$;
\medskip
Case $D_{23}$2.2. $-\frac13 \ge v$, $v\ge -\frac12u-\frac32$, $v\le -2u-3$;
\medskip
Case $D_{23}$2.3. $v\ge -\frac43$, $v\le -\frac12u-\frac32$, $v\ge -2u-3$;
\medskip
Case $D_{23}$2.4. $v\ge -\frac12u-\frac32$, $v\ge -2u-3$, $v\le -u-\frac53$;
\medskip
Case $D_{23}$2.5. $v\le -\frac12u-1$, $v\le -2u-2$, $v\ge -u-\frac53$
\medskip
The vertices of the Newton polygon for these cases are:
\settabs 6 \columns
\medskip
\+Case $D_{23}$2.1&$(1,-v)$& & $(4, 2u)$& &$(16, -4v-12)$\cr
\medskip
\+Case $D_{23}$2.2&$(1,-v)$& & $(4, 2u)$&$(8, -4v-8)$& $(16, -4v-12)$\cr
\medskip
\+Case $D_{23}$2.3&$(1,-v)$& $(2, -v-1)$& $(4, 2u)$&& $(16, -4v-12)$\cr
\medskip
\+Case $D_{23}$2.4&$(1,-v)$& $(2, -v-1)$& $(4, 2u)$&$(8, -4v-8)$& $(16, -4v-12)$\cr
\medskip
\+Case $D_{23}$2.5&$(1,-v)$& $(2, -v-1)$&& $(8, -4v-8)$& $(16, -4v-12)$\cr
\medskip
The ord's of $x_{j0}$ for these cases are:
\settabs 5 \columns
\medskip
\+Case $D_{23}$2.1&$\frac16u+\frac13v+1$&$\frac16u+\frac13v+1$& $-\frac23u-\frac13v$&$-\frac23u-\frac13v$ \cr
\medskip
\+Case $D_{23}$2.2&$\frac12$&$\frac12u+v+2$& $-\frac23u-\frac13v$&$-\frac23u-\frac13v$ \cr
\medskip
\+Case $D_{23}$2.3&$\frac16u+\frac13v+1$&$\frac16u+\frac13v+1$& $-u-\frac12v-\frac12$&1\cr
\medskip
\+Case $D_{23}$2.4&$\frac12$&$\frac12u+v+2$&$-u-\frac12v-\frac12$&1\cr
\medskip
\+Case $D_{23}$2.5&$\frac12$&$\frac12v+\frac76$&$\frac12v+\frac76$&1\cr
\medskip
We can rewrite the above information in notations $t=\ord a_2$ (numeration of cases is slightly another, sorry):
\medskip
Case $D_{23}$2.1. $-\frac43<v<-\frac13$, $t\le -v-3$ if $-\frac43<v\le-1$, and $t\le -4v-6$ if $-1\le v<-\frac13$;
\medskip
Case $D_{23}$2a.2. $-\frac43<v\le-1$, $-v-3 < t\le -4 v -6 $.
\medskip
Case $D_{23}$2a.3. $-\frac43<v<-1$, $-4v-6< t< -2 v -\frac{10}3$.
\medskip
Case $D_{23}$2b.2. $-1<v<-\frac13$, $-4v-6 < t< -v -3 $.
\medskip
Case $D_{23}$2b.3. $-1\le v<-\frac13$, $-v-3<t\le -2 v -\frac{10}3$.
\medskip
Case $D_{23}$2.5. $-\frac43<v<-\frac13$, $t\ge -2 v -\frac{10}3$.
\medskip
The vertices of the Newton polygon (for cases 2a, 2b):
\settabs 6 \columns
\medskip
\+Case $D_{23}$2.1&$(1,-v)$& & $(4, t)$& &$(16, -4v-12)$\cr
\medskip
\+Case $D_{23}$2a.2&$(1,-v)$& $(2,-v-1)$& $(4, t)$&& $(16, -4v-12)$\cr
\medskip
\+Case $D_{23}$2b.2&$(1,-v)$& & $(4, t)$&$(8, -4v-8)$& $(16, -4v-12)$\cr
\medskip
\+Case $D_{23}$2.3&$(1,-v)$& $(2,-v-1)$& $(4, t)$&$(8, -4v-8)$& $(16, -4v-12)$\cr
\medskip
\+Case $D_{23}$2.5&$(1,-v)$& $(2,-v-1)$& &$(8, -4v-8)$& $(16, -4v-12)$\cr
\medskip
Here we start a calculation of ord $x_1$ for some cases.
\medskip
Case $D_{23}$2.1, $j=1$. We have $\ord x_{10}=\frac16u+\frac13v+1$,
\medskip
$\ord b_{14}x_{10}^{16}=\frac83u+\frac43v+8$,
\medskip
$\ord b_{13}x_{10}^{8}=\frac43u-\frac43v+4$,
\medskip
$\ord b_{12}x_{10}^{4}=\frac23u+\frac13v+4$.
\medskip
In all points of $D_{23}$2.1 we have $\frac83u+\frac43v+8<\frac43u-\frac43v+4$,
hence, the domain $D_{23}$2.1 consists of two subdomains $D_{23}$2.1.1 ($v<-2u-4$) and $D_{23}$2.1.2 ($v>-2u-4$). We have the ord of the tail is $\frac83u+\frac43v+8$ in $D_{23}$2.1.1, $\frac23u+\frac13v+4$ in $D_{23}$2.1.2.
\medskip
Case $D_{23}$2.1.1. The line defined by points $(4, 2u)$, $(16, -4v-12)$ (vertices of Newton polygon) crosses the $v$-axis at the point $(0, \frac83u+\frac43v+4)$. If $u\ll0$ then the point $(1,-v)$ is over the segment $[(0, \frac83u+\frac43v+8), (4,2u)]$, hence this segment is a part of the Newton polygon of $x_1$, and $\ord x_{11}=\frac16u+\frac13v+2$. Let us continue for $\ord x_{12}$:
\medskip
$\ord b_{14}x_{11}^{16}=\frac83u+\frac43v+24$,
\medskip
$\ord b_{13}x_{11}^{8}=\frac43u-\frac43v+12$,
\medskip
$\ord b_{12}x_{11}^{4}=\frac23u+\frac13v+8$,
\medskip
$\ord b_{24}x_{10}^{16}=\frac83u+\frac43v+16$.

Apparently (to check ! ) continuing we get that domain $D_{23}$2.1.1 is uniformizable.
\medskip
Case $D_{23}$2.1.2. We have $\frac23u+\frac13v+4>\frac83u+\frac43v+4$ in this domain, hence if the point $(1,-v)$ is over the segment $[(0, \frac23u+\frac13v+4), (4,2u)]$ then $\ord x_{11}=-\frac13u+\frac1{12}v+1$. In any case, $\ord x_{11}\ge-\frac13u+\frac1{12}v+1$. Again apparently (to check ! ) continuing we get that domain $D_{23}$2.1.2 is uniformizable.
\medskip
For the cases $D_{23}$2a.3, $D_{23}$2b.2, $D_{23}$2b.3, $D_{23}$2.5 we have $h^1=4$, see Proposition 8.3.
\medskip
{\bf Case }$D_{23}$3. If $v\ge -\frac13$ the convex hull of $(1,-v)$, $(2, -v-1)$, $(8, -4v-8)$, $(16, -4v-12)$ consists of points $(1,-v)$, $(8, -4v-8)$, $(16, -4v-12)$. Depending on $u$ we have 2 domains on $(u,v)$-plane defined as below. Further, if we let $t=\ord a_2$ we get the third domain:
\medskip
Case $D_{23}$3.1. $v\ge -\frac13$, $v\le -\frac12u-\frac32$;
\medskip
Case $D_{23}$3.2. $v\ge-\frac13$, $-\frac12u-1\ge v\ge -\frac12u-\frac32$.
\medskip
Case $D_{23}$3.3. $v\ge-\frac13$, $t\ge -\frac{16}7 v -\frac{24}7$.
\medskip
Vertices of Newton polygon:
\settabs 6 \columns
\medskip
\+Case $D_{23}$3.1&$(1,-v)$& & $(4, 2u)$& &$(16, -4v-12)$\cr
\medskip
\+Case $D_{23}$3.2&$(1,-v)$& & $(4, 2u)$&$(8, -4v-8)$& $(16, -4v-12)$\cr
\medskip
\+Case $D_{23}$3.3&$(1,-v)$& & &$(8, -4v-8)$& $(16, -4v-12)$\cr
\medskip
Ord's of $x_{j0}$:
\settabs 5 \columns
\medskip
\+Case $D_{23}$3.1&$\frac16u+\frac13v+1$&$\frac16u+\frac13v+1$& $-\frac23u-\frac13v$&$-\frac23u-\frac13v$ \cr
\medskip
\+Case $D_{23}$3.2&$\frac12$&$\frac12u+v+2$& $-\frac23u-\frac13v$&$-\frac23u-\frac13v$ \cr
\medskip
\+Case $D_{23}$3.3&$\frac12$&$\frac37v+\frac87$&$\frac37v+\frac87$&$\frac37v+\frac87$ \cr
\medskip
The subdomain $D_{23}$3.3 is treated in Proposition 8.3.
\medskip

{\bf References}

\medskip
[A] Anderson Greg W. $t$-motives. Duke Math. J. Volume 53, Number 2 (1986), 457-502.
\medskip
[Dr]  Drinfeld, V. G. Elliptic modules. Math. USSR-Sb. no. 4, 561 - 592 (1976).
\medskip
[G] Goss, David Basic structures of function
field arithmetic.
Springer-Verlag, Berlin, 1996. xiv+422 pp.
\medskip
[GL07] Grishkov, A., Logachev, D. Duality of Anderson t-motives. arxiv.org/pdf/0711.1928.pdf
\medskip
[GL17] Grishkov, A., Logachev, D. Lattice map for Anderson t-motives: first approach. J. of Number Theory. 2017, vol. 180, p. 373 -- 402.
https://arxiv.org/pdf/1109.0679.pdf
\medskip
[GL18] Grishkov, A., Logachev, D. $h^1\ne h_1$ for Anderson t-motives. 2018. https://arxiv.org/pdf/1807.08675.pdf
\medskip
[H] Urs Hartl, Uniformizing the Stacks of Abelian Sheaves.

http://arxiv.org/abs/math.NT/0409341
\medskip
[L] Logachev, D. Anderson t-motives are analogs of abelian varieties with multiplication by imaginary quadratic fields. arxiv.org/pdf/0907.4712.pdf
\enddocument